\documentclass[11pt,a4text]{article}

\usepackage[latin1]{inputenc}
\usepackage[T1]{fontenc}
\usepackage[english]{babel}
\usepackage{ae}
\usepackage{amsthm, amssymb, amsmath, amsfonts, amscd}
\usepackage{graphicx,url}
\usepackage{bbm}

\newtheorem{theorem}{Theorem}

\newtheorem{lemma}{Lemma}
\newtheorem{proposition}{Proposition}
\newtheorem{corollary}{Corollary}

\newtheorem{remark}{Remark}

\newcommand{\Hd}{H_\delta}
\newcommand{\Hdn}{H_{\delta_n}}
\newcommand{\Hdk}{H_{\delta_k}}
\newcommand{\Hud}{H_{\delta}^{u}}
\newcommand{\Hdd}{H_{\delta}^{d}}
\newcommand{\Hld}{H_{\delta}^{l}}
\newcommand{\Hrd}{H_{\delta}^{r}}
\newcommand{\Hudp}{\tilde{H}_{\delta}^{u}}
\newcommand{\Hddp}{\tilde{H}_{\delta}^{d}}
\newcommand{\Hldp}{\tilde{H}_{\delta}^{l}}
\newcommand{\Hrdp}{\tilde{H}_{\delta}^{r}}
\newcommand{\Hudk}{H_{\delta_k}^u}
\newcommand{\Hddk}{H_{\delta_k}^d}
\newcommand{\Hldk}{H_{\delta_k}^l}
\newcommand{\Hrdk}{H_{\delta_k}^r}
\newcommand{\Nld}{N_\delta^l}
\newcommand{\Nrd}{N_\delta^r}
\newcommand{\Expect}{\mathbb{E}}
\newcommand{\Qud}{Q_\delta^u}
\newcommand{\Qdd}{Q_\delta^d}
\newcommand{\onedge}{\vec{e}}
\newcommand{\oes}{\vec{e^*}}
\newcommand{\jes}{\vec{(\tau.e)^*}}
\newcommand{\jses}{\vec{(\tau^2.e)^*}}
\newcommand{\Prob}{\mathbb{P}}
\newcommand{\jedge}{\vec{\tau . e}}
\newcommand{\jsedge}{\vec{\tau^2 . e}}
\newcommand{\Omdel}{\Omega_\delta}
\newcommand{\diffoe}{\partial_{\onedge}}
\newcommand{\diffoep}{\partial_{\onedge}^+}
\newcommand{\diffoepm}{\partial_{\onedge}^{\pm}}
\newcommand{\diffmoep}{\partial_{-\onedge}^+}
\newcommand{\diffoem}{\partial_{\onedge}^-}

\newcommand{\diffjep}{\partial_{\jedge}^+}

\newcommand{\diffjsep}{\partial_{\jsedge}^+}
\newcommand{\diffjsem}{\partial_{\jsedge}^-}
\newcommand{\LHex}{\mathcal{L}}
\newcommand{\RHex}{\mathcal{R}}
\newcommand{\IHex}{\mathcal{I}}
\newcommand{\THex}{\mathcal{T}}
\newcommand{\DefHex}{\mathcal{H}}
\newcommand{\uside}{u}
\newcommand{\dside}{d}

\newcommand{\Lww}{\LHex_w^w}
\newcommand{\Rbb}{\RHex_b^b}
\newcommand{\Rww}{\RHex_w^w}
\newcommand{\Ibb}{\IHex_b^b}
\newcommand{\Iww}{\IHex_w^w}
\newcommand{\dz}{\, \mathrm{d}z}
\newcommand{\gamdel}{\gamma_{\delta}}
\newcommand{\gamint}{\mathrm{Int}(\gamma_{\delta})}
\newcommand{\gamhex}{\mathrm{Hex}(\gamma_\delta)}
\newcommand{\hexboundary}{\partial f}
\newcommand{\hcst}{\frac{\sqrt{3}}{2} i}

\title{Critical percolation: the expected number of clusters in a rectangle}
\author{Cl\'{e}ment Hongler and Stanislav Smirnov \\ Universit\'{e} de Gen\`{e}ve}

\begin{document}
\maketitle

\begin{abstract}
We show that for critical site percolation on the triangular lattice
two new observables have conformally invariant scaling limits.
In particular the expected number of clusters separating two pairs of points
converges to an explicit conformal invariant.
Our proof is independent of earlier results and $SLE$ techniques,
and in principle should provide a new approach to establishing conformal
invariance of percolation.
\end{abstract}

\section{Introduction}

Percolation is perhaps the easiest  two-dimensional lattice model to formulate,
yet it exhibits a very complicated behavior.
A number of spectacular predictions (unrigorous, but very convincing)
appeared in the physics literature over the last few decades, see \cite{Cardy}.
One of them, the Cardy's formula for the scaling limit of crossing probabilities,
was recently established for the critical site percolation on triangular lattice
\cite{Smirnov}.
Consequently, scaling limits of interfaces were identified with Schramm's
$SLE_6$ curves, and many other predictions were proved, see e.g. \cite{Smirnov-Werner}.

In this paper we show that two new observables 
for the critical site percolation on triangular lattice
have conformally invariant scaling limits.
Furthermore, we obtain explicit formulae, consistent with predictions obtained by physicists
\cite{Cardy-expected,Simmons}.
Our proof is independent of earlier conformal invariance results,
and uses methods similar to those in \cite{Smirnov} rather than $SLE$ techniques.
It is also restricted to the same triangular lattice model,
but at least one should be able to use it for a new proof
of the conformal invariance in this case.

\subsection{Acknowledgements}
The first author would like to thank Thomas Mountford
and Yvan Velenik for useful discussions and remarks. 
This work was supported by the Swiss National Science Foundation grants 117596, 117641, 121675. 
The first author was partially supported by an EPFL Excellence Scholarship.

\section{Notation and Setup}

For convenience reasons, in this paper we shall not work on the triangular lattice but rather on its dual,
the \emph{honeycomb lattice}, and thus, rather than coloring vertices of triangles, we shall color hexagonal faces
(which is obviously equivalent).

\subsection{Graph and Model}

Let $ \Omega \subset \mathbb{C} $ be a \emph{Jordan domain} (whose boundary is a simple closed curve), and orient its
boundary $ \partial \Omega $ counterclockwise. Let $ l $ and $ r $ be two distinct points on $ \partial \Omega $,
which separate it into a curve $ \uside $ going from $ r $ to $ l $ (with respect to the orientation
of $ \partial \Omega $) and a curve $ \dside $ going from $ l $ to $ r $, such that we have
$ \partial \Omega \setminus \{ l, r \} = \uside \cup \dside $. Let finally $ w $ be a point on $ \uside $.

\begin{remark}
	The assumption on $ \Omega $ to be a Jordan domain is not really necessary, and the result remains true
	under weaker assumptions detailed in Section 3.
	We use this assumption in Section \ref{precompactness} to avoid lengthy and not so interessant discussions.
\end{remark}

\begin{figure}[!ht]
\centering
\includegraphics[width=12cm]{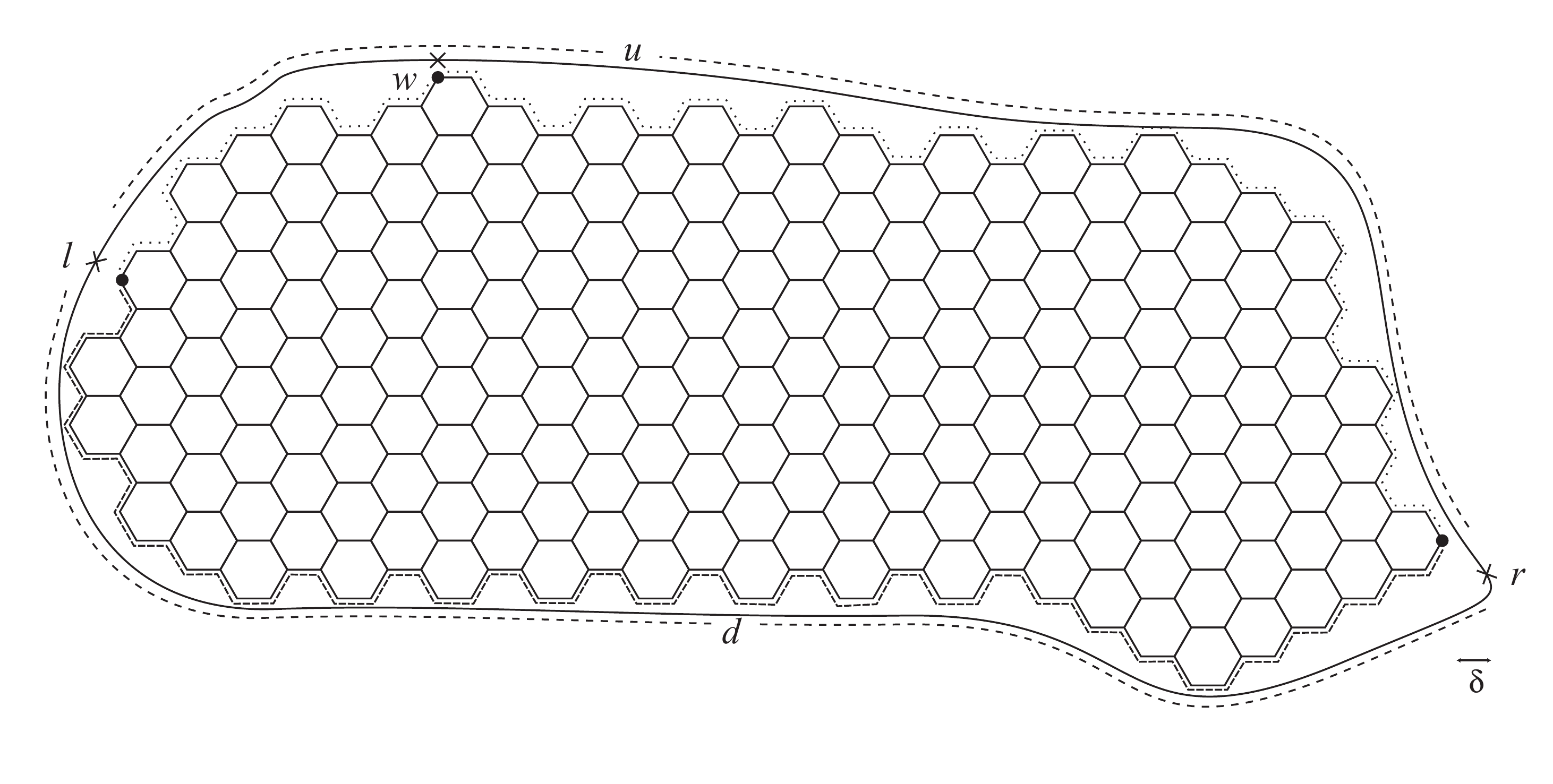}
\caption{Domain discretization: the points are marked with crosses, while their discretization are
marked with circles. The discretization of $ \uside $ and $ \dside $ are depicted by sparse and
dense dashes respectively.}
\end{figure}

We consider the \emph{discretization $ \Omdel $ of $ \Omega $} by regular hexagons defined as follows. 
Let $ G_\delta $ be the regular hexagonal lattice embedded in the complex plane with mesh size
(i.e. sidelength of a hexagon) $ \delta > 0 $.
We define $ \Omdel $ as the graph obtained by taking a maximal connected 
component made of hexagonal faces of $ G_\delta $: the union of the closure of the faces 
is a simply connected subset of $ \mathbb{C} $. We denote by $ \bar{\Omega}_{\delta} $ this subset
and by $ \partial \Omdel $ the (counterclockwise-oriented) 
simple path consisting in edges of $ \Omdel $ such that $ \Omdel $ is contained inside it. We define
the discretization of $ l $, $ r $ and $ w $ as the closest corresponding vertices of $ \partial \Omdel $, and
$ \uside $ and $ \dside $ as the paths from $ r $ to $ l $ and $ l $ to $ r $ respectively, following the
orientation of $ \partial \Omdel $. In general we will identify $ l, r, \uside, \dside $ to their
respective discretizations.

We are interested in the process of \emph{critical percolation on the faces of $ \Omdel $}: 
each face of $ \Omdel $ is
colored either in white or in black with probability $ \frac{1}{2} $, independently of the other faces; 
such a coloring is called a \emph{configuration}.
More precisely, we are interested in the \emph{scaling limit} of this model: the description of the
global geometry of the process as the mesh size $ \delta $ tends to $ 0 $.

Note that for this model $ \frac{1}{2} $ is known to be the critical value for the probability
thanks to the work of Wierman and Kesten.
However, we do not use that this value is critical, only that it is self-dual. 

We call \emph{path of hexagons} a sequence $ \DefHex_1, \ldots, \DefHex_n $ of hexagons such that
$ \DefHex_i $ is adjacent to $ \DefHex_{i + 1} $ for $ i = 1, \ldots, n - 1 $; a path is called
\emph{simple} if all of its hexagons are distinct; a \emph{closed} simple path (the last hexagon is adjacent
to the first one) is called a \emph{circuit}. 
We say that a hexagon $ \DefHex $ is \emph{connected to $ \uside $ by a white path} 
if there exists a path of white hexagons
that contains $ \DefHex $ and that hits $ \uside $ (contains a hexagon having an edge belonging
to (the discretization of) $ \uside $). We define
similarly connection events involving black instead of white paths or connections to $ \dside $ instead of $ \uside $.
We say that a path of hexagons $ \gamma $ \emph{separates} 
two families of points $ A $ and $ B $ if the interior of each continuous path
$ \alpha $
contained in $ \bar{\Omega}_{\delta} $ from a point of $ A $ to a point of $ B $ crosses 
the closure of a hexagon of $ \gamma $.

We call \emph{white cluster} a connected (i.e. path-connected in the sense defined above) set of white hexagons. 
For a cluster $ K $ touching $ \uside $ and $ \dside $, we define its \emph{left boundary}
(respectively \emph{right boundary}) as the left-most (respectively right-most) simple path contained in $ K $
that touches $ \uside $ and $ \dside $, i.e. such that there is no path in $ K $ separating it from $ l $ (respectively $ r $); elementary topological considerations show that this notion is well-defined.
One important property of our lattice is indeed its \emph{self-duality}: the boundary of a
white cluster (that does not touch the boundary) is a black circuit, and vice versa.

We will use the term left boundary for simplicity,
but strictly speaking our definition gives the left-most simple curve inside the cluster,
that is its left-most boundary after ``peninsulas'' attached by only one hexagon are erased.
So this curve would rather bound on the right the dual cluster bordering ours on the left.

\begin{figure}[!ht]
\centering
\includegraphics[width=13cm]{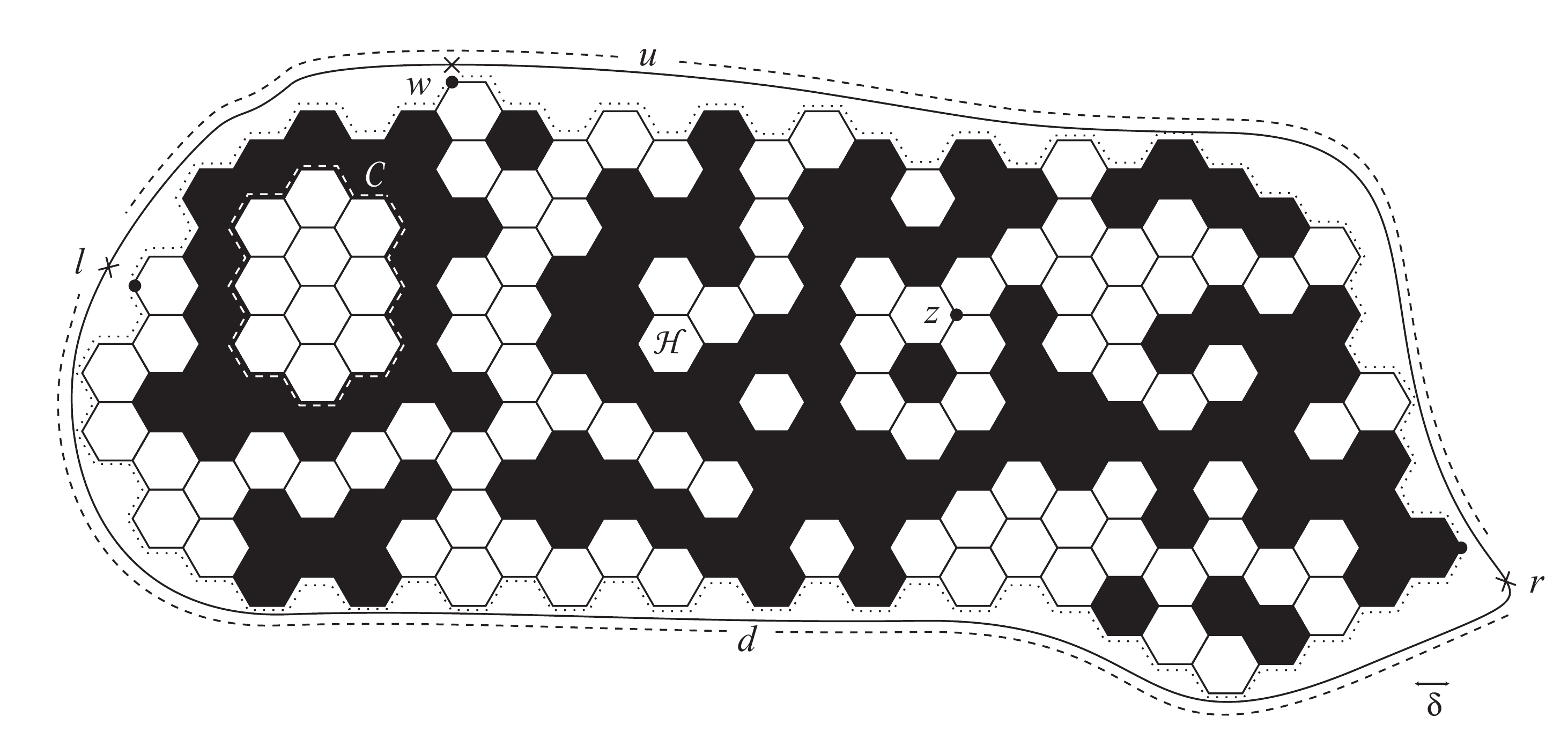}
\caption{In this configuration, the dashed region $ C $ is a white cluster, the hexagon $ \DefHex $ is
connected to $ \uside $ by a white path and the points $ l $ and $ w $ are separated from the points
$ z $ and $ r $ by a black path.}
\end{figure}

Notice that since the probability for a hexagon to be white is $ \frac{1}{2} $, 
any event (i.e. set of configurations) has the same probability as its negative with
respect to the colors: for instance, the probability that there is a white path from $ x $ to $ y $ is the same
as the one that there is a black path from $ x $ to $ y $. For an event $ A $, we will denote by $ \tilde{A} $ the
\emph{negative event}: a configuration $ \omega $ belongs to $ A $ if and only if the \emph{negative configuration}
$ \tilde{\omega} $ (i.e. with the colors black and white flipped) belongs to $ \tilde{A} $.

\subsection{Observables}

Let $ \delta > 0 $ and consider the process of percolation on $ \Omdel $ as described in the previous section.
For each $ z $ vertex of $ \Omdel $ we define the following random variables and events:
\begin{itemize}
	\item $ \Nld ( z ) $: the number of (simple) left boundaries of white clusters touching
	$ \uside $ and $ \dside $ separating $ l $ and $ w $ from $ z $ and
	$ r $ \emph{minus} the number of (simple) left boundaries of white clusters (touching
	$ \uside $ and $ \dside $) separating $ l $ and $ z $ from $ w $ and $ r $;
	\item $ \Nrd ( z ) $: the same as $ \Nld ( z ) $ but for (simple) right boundaries of 
	white clusters (also touching $ \uside $ and $ \dside $);
	\item $ \Qud ( z ) $: the event that there exists a white simple path from $ \dside $ to 
	$ \dside $ that separates $ z $ from $ l $ and $ r $
	and that is connected to $ \uside $;
	\item $ \Qdd ( z ) $: the same event as $ \Qud ( z ) $ but with a
	white simple path path from $ \uside $ to $ \uside $ connected to $ \dside $ instead.
\end{itemize}
This allows us to define our observables:
\begin{eqnarray*}
	\Hld ( z ) := \Expect [ \Nld ( z ) ], & \quad &
	\Hrd ( z ) := \Expect [ \Nrd ( z ) ], \\
	\Hud ( z ) := \Prob [ \Qud ( z ) ], & \quad &
	\Hdd ( z ) := \Prob [ \Qdd ( z ) ].
\end{eqnarray*}
We extend these functions to continuous functions on $ \bar{\Omega}_{\delta} $ in the following way (in fact any 
reasonable manner will work): first for the center of a hexagon, take the average value of its vertices. Then
divide the hexagon into six equilateral triangles, and define the functions on each triangle by affine interpolation.
We can then extend the functions to $ \bar{\Omega} $ in a smooth way.

\begin{remark}
	The point $ w $ could in fact be anywhere in $ \bar{\Omega} $ (changing its position only modifies
	the functions $ \Hld $ and $ \Hrd $ by an additive constant). 
	In our setup it lies on the boundary for simplicity.
\end{remark}
\begin{remark}
	
	Another way of writing $ \Hld $ (similarly for $ \Hrd $), which motivates its definition, is the following:
	count the expected number of left boundaries that separate $ l $ from $ z $ and $ r $ \emph{minus} the expected 
	number of left boundaries that separate $ l $ from $ w $ and $ r $. 
	
	It is easy to check that this definition is equivalent to the one given above (the boundaries that count
	positively are precisely the ones that separate $ l $ from $ r $ and $ z $ but not $ w $, the boundaries
	that count negatively are the ones that separate $ l $ from $ r $ and $ w $ but not $ z $).
	
	If one uses this way to write $ \Hld $, taking the difference is essential to get a finite limit: 
	as the mesh tends to $ 0 $ the expected number of clusters joining $ \uside $ to $ \dside $ blows up.
	
\end{remark}
\begin{remark}
	
	Notice that the quantities $ \Hld $ and $ \Hrd $ are the same: if one has a configuration in
	$ \{ \Nld ( z ) = k \} $, flipping the colors of all the hexagons gives a 
	configuration in $ \{ \Nrd ( z ) = k \} $.
	
\end{remark}

\begin{figure}[!ht]
\centering
\includegraphics[width=12cm]{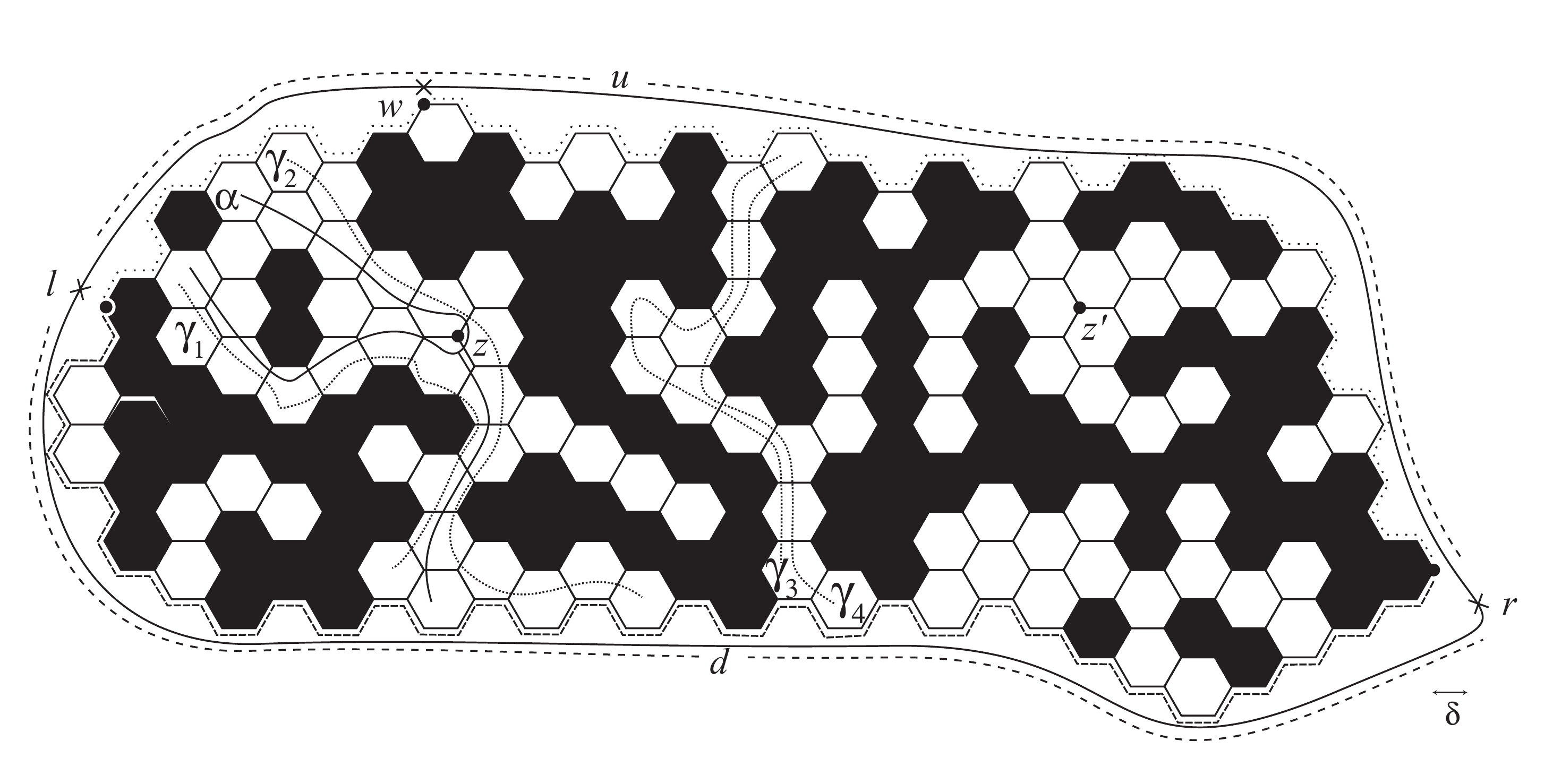}
\caption{In this configuration, $ \gamma_1 $ and $ \gamma_3 $ are left boundaries of white clusters,
$ \gamma_2 $ and $ \gamma_4 $ are right boundaries. $ \Qdd ( z ) $ occurs and we 
have $ \Nld ( z ) = 0, \Nrd ( z ) = -1 $ and $ \Nld ( z' ) = \Nrd ( z' ) =1 $.}
\end{figure}

\section{Conformal invariance and main result}

By conformal invariance of critical percolation we mean that the same observable on two conformally
equivalent Riemann surfaces has the same scaling limit.

It was proven in \cite{Smirnov} that crossing probabilities of conformal rectangles (here the Riemann
surface is a simply connected domain with four marked boundary points) are conformally invariant and
satisfy Cardy's prediction.

Consequently the interfaces of macroscopic clusters converge to Schramm's SLE curves and we can deduce
conformal invariance of many other observables.

The goal of this paper is to show conformal invariance of the observables
$ \Hld + \Hrd $ and $ \Hud - \Hdd $ in the same setup, without appealing to the results of \cite{Smirnov}.

\subsection{Limit of the observables}

In order to get our conformal invariance result, 
we prove a more geometrical one: a linear combination of our two observables
turns out to be (in the limit) a conformal mapping. 
For each $ \delta > 0 $, define $ \Hd := ( \Hld + \Hrd ) - \frac{\sqrt{3}}{2}i ( \Hud - \Hdd ) $. Then we have:

\begin{theorem} \label{mainthm}
	As $ \delta $ tends to $ 0 $, $ \Hd $ converges uniformly on the compact
	subsets of $ \bar{\Omega} \setminus \{ l , r \} $
	to a function $ h $ which is the unique conformal mapping from $ \Omega $ to the strip 
	$ S: = \{ x + iy : x \in \mathbb{R}, y \in ( -\sqrt{3} / 4, \sqrt{3} / 4) \} $ that maps
	(in the sense of prime ends) $ l $ to the left end of the strip ,
	$ r $ to the right end of the strip and $ w $ to $ \frac{\sqrt{3}}{4} i $.	
\end{theorem}

\begin{remark}
	The theorem remains valid under the weaker assumption that the discretizations
	$ \Omdel $ of the domain converge in Caratheodory's sense to $ \Omega $, in which case
	the observables converge on the compact subsets of $ \Omega $.
\end{remark}

This theorem gives us the asymptotical conformal invariance (and the existence of the limit)
of the two observables $ \Hld + \Hrd $ and $ \Hud - \Hdd $ in the following sense. 

\begin{corollary} \label{confinv}
	Let $ \Phi $ be a conformal map as above and denote by $ \Hldp $,  $ \Hrdp $, $ \Hudp $ and $ \Hddp $ 
	the corresponding observables on
	the domain $ \Omega' := \Phi ( \Omega ) $ 
	with the corresponding points $ l' := \Phi(l), r' := \Phi(r), w':= \Phi(w) $. Then we have the following
	conformal invariance result:
	\begin{eqnarray*}
		\lim_{\delta \to 0} \Hld + \Hrd & = & \lim_{\delta \to 0} (\Hldp + \Hrdp) \circ \Phi \\
		\lim_{\delta \to 0} \Hud - \Hdd & = & \lim_{\delta \to 0} (\Hudp - \Hddp) \circ \Phi
	\end{eqnarray*}
\end{corollary}

\begin{proof}
	By uniqueness of the conformal mapping to $ S $ with three points fixed we have
	$ h = h ' \circ \Phi $ (the images of $ l $, $ r $ and $ w $ by $ h $ and $ h ' \circ \Phi $ are the same). 
	Taking the real and imaginary parts gives the result.
\end{proof}

Taking $ z $ and $ w $ on the boundary 
we obtain the conformal invariance of the expected number of clusters in a \emph{conformal rectangle}
(a Jordan domain with four distinct points on its boundary).
Let $ \Xi $ be a conformal rectangle with the four points $ a_1, a_2, a_3, a_4 $ 
in counterclockwise order.
Discretize the domain and the four points as before and consider the expected number
$ C_\delta $ of white clusters separating
$ a_1 $ and $ a_4 $ from $ a_2 $ and $ a_3 $, counted in the following way:
\begin{itemize}
	\item If a cluster touches both (the discretization of) the arcs $ a_4 a_1 $ and $ a_2 a_3 $ (along the counterclockwise
	orientation of the $ \partial \Xi $), it does not count.
	\item If a cluster touches exactly one of the arcs $ a_4 a_1 $ and $ a_2 a_3 $, it counts once.
	\item If it does not touch any of the two arcs, it counts twice.
\end{itemize}

\begin{corollary}\label{clustnumb}

	The quantity $ C_\delta $ admits a conformally invariant limit as $ \delta \to 0 $:
	If $ \Xi' $ is another conformal rectangle with the four points $ a_1', a_2', a_3', a_4' $, if 
	$ \Psi : \Xi \to \Xi' $ is a conformal mapping such that $ \Psi (a_i) = a_i' $ for $ i = 1, 2, 3, 4 $,
	and $ C_\delta ' $ is the corresponding number in the domain $ \Xi ' $, then
	\[
		\lim_{\delta \to 0} C_\delta = \lim_{\delta \to 0} C_\delta '
	\]
	
\end{corollary}
\begin{proof}
	It suffices to take $ z $ on the boundary (choose
	$ z = a_1, w = a_2, l = a_3, r = a_4 $) and to see that in this case $ C_\delta = \Hld + \Hrd $:
	no clusters count negatively, 
	if a cluster does not touch any arc, both its left and right boundaries count, etc.
	Therefore the result follows from the previous corollary.
\end{proof}

\subsection{Formulae}

It is not difficult to express the quantity $ C_\delta $ in terms of the cross-ratio (the conformal map
from the half-plane to a strip is simply a logarithm). If we denote by 
$ \lambda = \frac{(a_1 - a_3)(a_2 - a_4)}{(a_1 - a_4)(a_2 - a_3)} $ the cross-ratio of the four points,
we get
\[
	\lim_{\delta \to 0} C_\delta = \frac{\sqrt{3} \pi}{2} \log \left(\frac{1}{1 - \lambda} \right).
\]

By adding to this formula
the probability that a cluster (separating $ a_4 $ and $ a_1 $ from $ a_2 $ and $ a_3 $) 
touches the arc $ a_4 a_1 $ and the 
probability that such a cluster touches moreover the arc $ a_2 a_3 $ one can obtain
(twice) the expected number of clusters without specific counting.

Using self-duality one can show that these two quantities are the same and that they can
be expressed as the difference of the probability that a cluster separates
$ a_4 $ and $ a_1 $ from $ a_2 $ and $ a_3 $
minus one half times the probability that a cluster touches the four sides of our conformal rectangle:
if there is a black cluster separating $ a_4 $ and $ a_1 $ from $ a_2 $ and $ a_3 $ (event $ \mathcal{E}_1 $ on
Figure 5), then consider the right-most such cluster;
either it touches also the arc $ a_1 a_4 $ 
(event $ \mathcal{E}_3 $) or it does not; in this latter case by self-duality
there is a white cluster on its right touching the arcs $ a_3 a_4 $, $ a_4 a_1 $ and $ a_1 a_2 $ 
(event $ \mathcal{E}_2 $). 

Then we can decompose the event $ \mathcal{E}_3 $ in the following way.
Either the cluster touching the arcs $ a_3 a_4 $, $ a_4 a_1 $ and $ a_1 a_2 $ touches also
arc $ a_2 a_3 $ (event $ \mathcal{E}_5 $) or it does not and there is a white cluster that separates
it from arc $ a_2 a_3 $ (event $ \mathcal{E}_4 $). 

A color-flipping argument gives that
events $ \mathcal{E}_2 $ and $ \mathcal{E}_4 $ have the same probability (one has that the negative
$ \tilde{\mathcal{E}}_2 $ of $ \mathcal{E}_2 $ is $ \mathcal{E}_4 $), which is therefore
$ \frac{1}{2} \left( \Prob [ \mathcal{E}_1 ] - \Prob [ \mathcal{E}_5 ] \right) $. Since by
self-duality $ \Prob [ \mathcal{E}_3 ] = \Prob [ \mathcal{E}_4 ] + \Prob [ \mathcal{E}_5 ] $ we
obtain $ \Prob [ \mathcal{E}_3 ] = \frac{1}{2} \left( \Prob [ \mathcal{E}_1 ] + \Prob [ \mathcal{E}_5 ] \right) $.

Both quantities are conformally invariant and given by Cardy's formula (see \cite{Cardy}, and
\cite{Smirnov} for a proof) and by Watts'
formula respectively (see \cite{Watts}, and \cite{Dubedat} for a proof).

\begin{figure}[!ht]
\centering
\includegraphics[width=12cm]{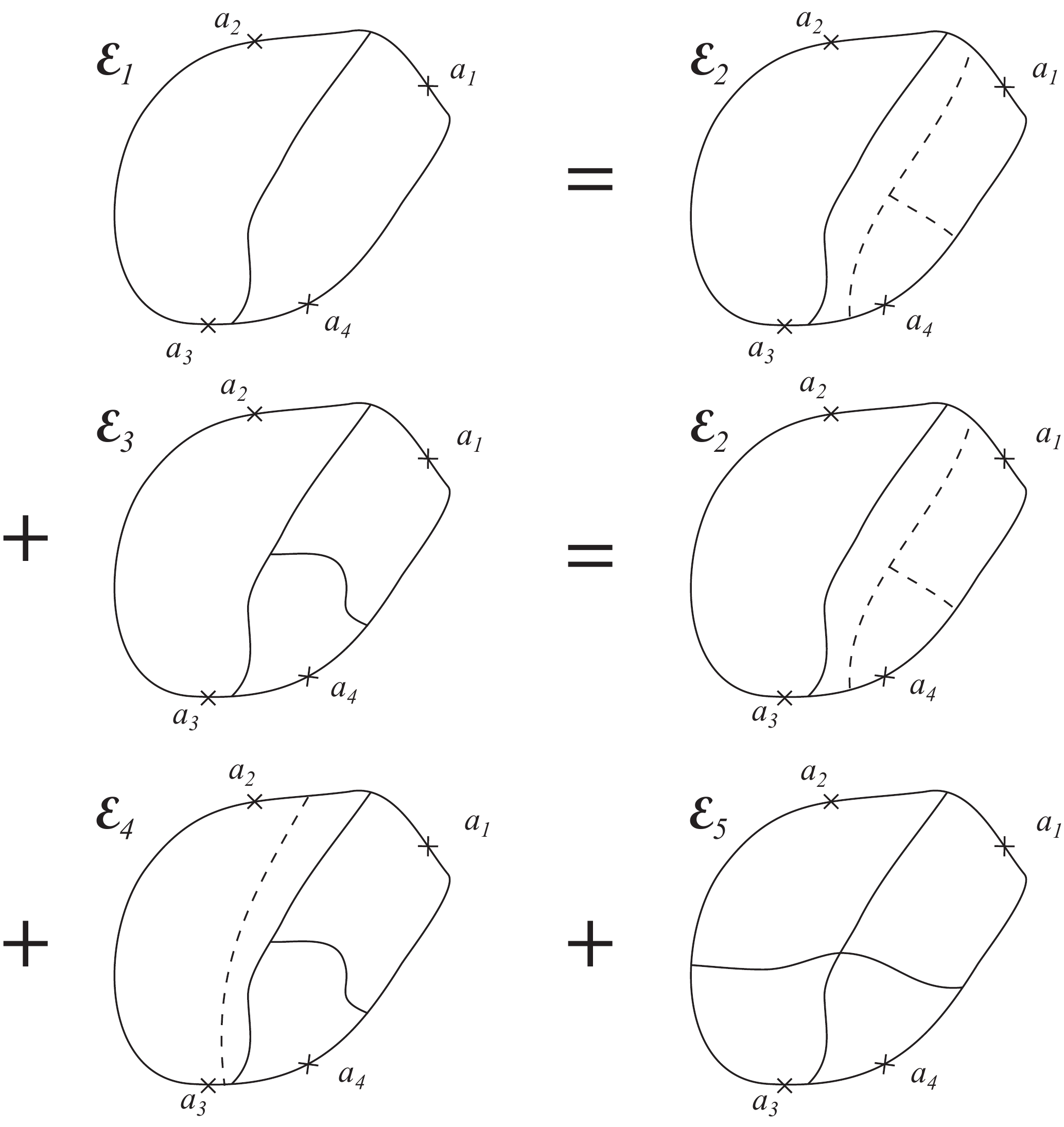}
\caption{Self-duality decomposition}
\end{figure}

So one obtains eventually:

\begin{proposition} \label{exp-num-clust}
The scaling limit of the the expected number of clusters separating
$ a_1 $ and $ a_4 $ from $ a_2 $ and $ a_3 $ is equal to:
\[
	\frac{2 \pi \sqrt{3}}{\Gamma (\frac{1}{3})^3} \lambda^{\frac{1}{3}}
	\, \, {_2 F _1} \left( \frac{1}{3}, \frac{2}{3} ; \frac{4}{3} ; \lambda \right)
	- \frac{1}{2} \frac{\sqrt{3}}{2 \pi}  \lambda \, \, {_3 F _2} 
	\left( 1, 1, \frac{4}{3} ; \frac{5}{3}, 2; \lambda \right)
	+ \frac{\sqrt{3}}{4 \pi} \log \left( \frac{1}{1 - \lambda} \right),
\]
where the first term comes from Cardy's formula, the second from Watts' formula and the third from
the main result of our paper.
\end{proposition}

\subsection{Open questions}
In this paper we show that certain observables have conformally invariant scaling limits.
The most prominent mathematical tool for rigorous treatment of conformal invariance 
is Schramm's $SLE$, which describes scaling limits of interfaces 
by the traces of the randomly driven Loewner evolutions --
the so called $ SLE_{\kappa} $ curves,
see \cite{Lawler} for an introduction.
Once convergence to $ SLE $ is known, many quantities related to the model can be computed.
The only proof for percolation uses Cardy's formula for crossing probabilities
(established for triangular lattice only in \cite{Smirnov}) 
and locality of percolation or the so-called ``martingale trick'',
see \cite{Smirnov-preprint,Smirnov-icmp,Camia}.

\subsubsection{How to use our observables to establish conformal invariance  of  critical percolation?}

Whether our observable can replace the crossing probabilities in the proofs above,
is interesting even if it has no less dependence on the triangular lattice. 
The problem that prevents the direct application  of the same technique as in \cite{Smirnov-icmp}
is that our observable does not have a ``martingale'' property 
(see \cite{Smirnov2} for an overview)
with respect to the percolation interface.
However, one can attempt other approaches, for example exploiting locality.

\subsubsection{Are our observables computable with $ SLE $?}

For the same reason, computing our observables
with $ SLE $ techniques
(using this time that the percolation scaling limit is described by $ SLE_6 $)
is not immediate. 
In principle, the computation should be possible, but the setup might be difficult. 

\subsubsection{Are there other similar observables?}

Similar techniques allow to compute crossing probabilities and two similar observables
in this paper.
One can ask how much more one can learn without appealing to $SLE$ techniques,
in particular whether there are any other computable observables?

\section{Outline of the proof}

The proof of Theorem \ref{mainthm} consists of three parts. 

\begin{itemize}
	\item
	First we prove that from each sequence $ \left( \Hdn \right)_{n \in \mathbb{N}} $, with
	$ \delta_n $ tending to $ 0 $, one can extract a subsequence $ \Hdk $ which converges uniformly on the compact subsets
	of $ \bar{\Omega} \setminus \{ l , r \} $ to a limit function $ h $.
	
	\item 
	We show then that any such subsequential limit $ h $ satisfies the following boundary conditions:
	\begin{eqnarray*}
		\mathrm{Im} ( h ) & = & \frac{\sqrt{3}}{4} \, \, \, \mbox{ on }  \uside \\
		\mathrm{Im} ( h ) & = & - \frac{\sqrt{3}}{4} \, \, \, \mbox{ on }  \dside \\
		\mathrm{Re} ( h ( w ) ) & = & 0 
	\end{eqnarray*}
	
	\item 
	We prove finally that $ h $ is analytic.
\end{itemize}

In order to get that $ h $ is the conformal map $ \phi $ of Theorem \ref{mainthm}, we observe that $ h $ and
$ \phi $ have the same imaginary part (on the boundary and hence inside since the imaginary part is
harmonic), and thus have the same real part up to a (real) constant by Cauchy-Riemann
equations. The constant is $ 0 $ since
the real part of both is equal to $ 0 $ at $ w $. Since any subsequential limit has the desired value, we
conclude by precompactness that $ \Hd $ converges to $ \phi $.
\section{Precompactness}\label{precompactness}

In order to prove the precompactness of the family of functions $ ( \Hd )_{\delta > 0} $, we show that the four families
$ ( \Hld )_{\delta}, ( \Hrd )_{\delta}, ( \Hud )_{\delta}, ( \Hdd )_{\delta} $ are uniformly Hölder
continuous on each compact
subset of $ \bar{\Omega} \setminus \{ l, r \} $. Notice that since the interpolation is regular enough
we may suppose in the estimates that the points we are considering are vertices of the hexagonal faces.

\begin{lemma}\label{hudreg}
	For every compact $ K \subset \bar{\Omega} \setminus \{ l, r \} $, 
	the functions $ \Hud $ and $ \Hdd $ are uniformly H\"{o}lder continuous on $ K $ with respect to the metric
	$ d $ of the length of the shortest path in $ \bar{\Omega} $.
\end{lemma}

\begin{proof}
	We prove the result for $ \Hud $. 
	
	Let $ \beta = \inf_{z \in K} \left( \max \left( \mathrm{dist} (z, \uside), 
	\mathrm{dist} (z, \dside) \right) \right) $. By compactness of $ K $ we have $ \beta > 0 $, and
	so each point in $ K $ is at distance at least $ \beta $ from $ \uside $ or $ \dside $.
	
	We have that for each $ z, z' \in \Omega $ the disc $ D:= D((z+z')/2, d(z,z')) $ contains a path from
	$ z $ to $ z' $. 
	
	Since $ | \Hud | $ is uniformly bounded, we can assume from now that the points
	$ z $ and $ z ' $ (in $ K $) are close enough, i.e. such that
	$ d( z, z ') \leq \beta $.
	By elementary partitioning we have that 
	$ | \Hud ( z ) - \Hud ( z ') |
	\leq \Prob [ \Qud ( z ) \setminus \Qud ( z ') ] + \Prob [ \Qud ( z ' ) \setminus \Qud ( z) ] $. 
	So it is enough to show that there exists $ C > 0 $ and $ \alpha > 0 $ such that
	\[
		\Prob [ \Qud ( z ) \setminus \Qud ( z ' ) ] \leq C \cdot d (z, z')^{\alpha}.
	\]
	
	\begin{figure}[!ht]
\centering
\includegraphics[width=12cm]{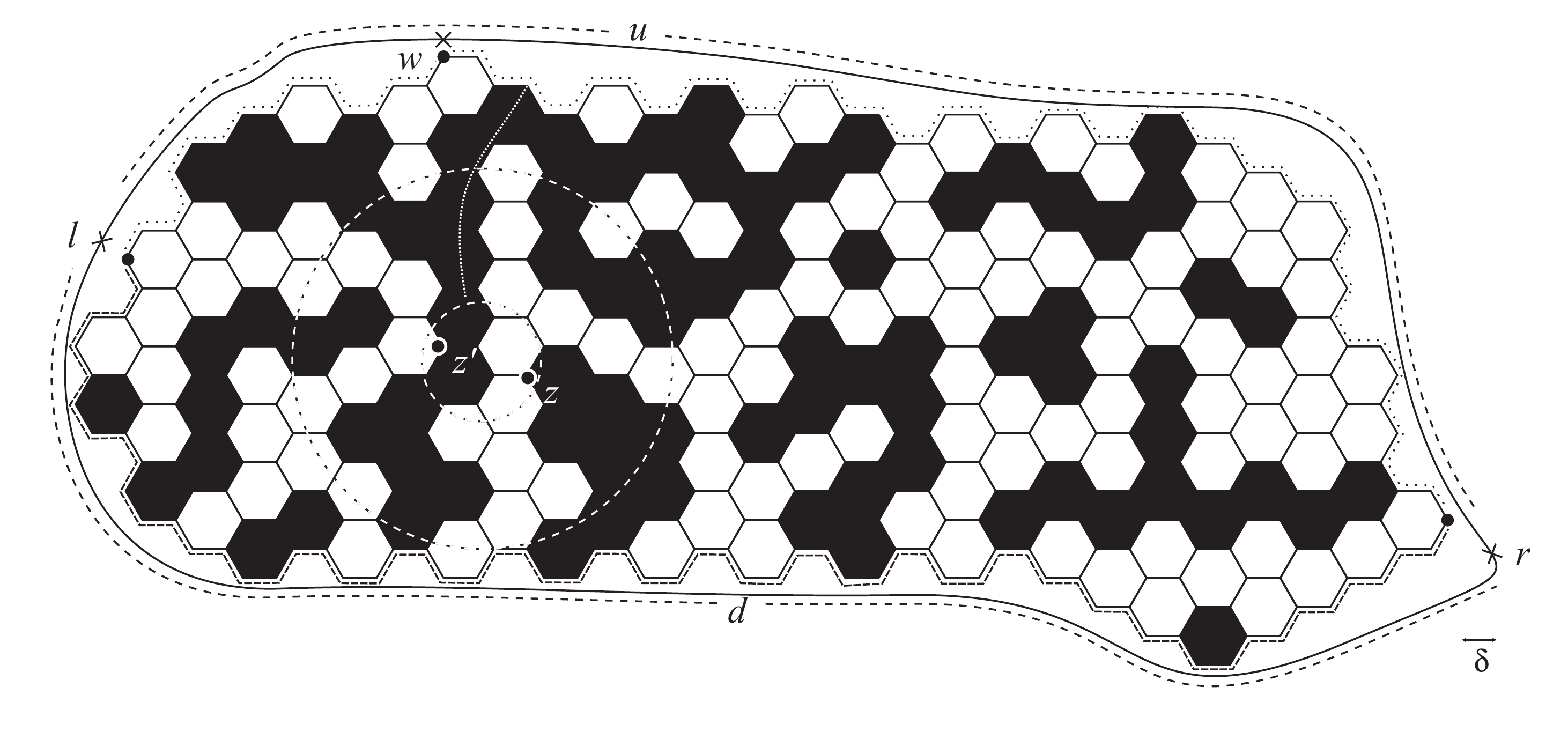}
\caption{The event $ \Qud ( z ) \setminus \Qud ( z' ) $ implies a black connection of a microscopic
circle of radius $ \ell ( | z - z' | ) $ to a macroscopic circle of radius $ \beta $ ($ K $ is not fixed precisely 
on this picture).}
\end{figure}

	By self-duality, we have that the occurence of the event $ \Qud ( z ) \setminus \Qud ( z ' ) $ implies the
	connection of the boundary of the disc $ D $ to $ \uside $ by a black path and to $ \dside $ by two disjoint white
	paths. Since at least one of the
	two sides is at distance $ \beta $ (for $ \delta $ sufficiently small, which we may suppose), this event implies
	the connection (by a black or white path) 
	of a (\emph{microscopic}) circle of radius $ d ( z, z') $ to a circle of
	(\emph{macroscopic}) radius $ \beta $. 
	By Russo-Seymour-Welsh Theorem (see \cite{Bollobas}, \cite{Grimmett} for instance), 
	there exists $ C > 0 $ and $ \alpha > 0 $
	such that this event is of probability less that $ C \cdot d(z, z')^{\alpha} $ (uniformly in $ \delta $)
	and this gives us the desired result.
\end{proof}

\begin{lemma}\label{hldreg}
	For every compact $ K \subset \bar{\Omega} \setminus \{ l, r \} $, 
	the functions $ \Hld $ and $ \Hrd $ are uniformly bounded and uniformly H\"{o}lder continuous on $ K $ with respect
	the metric $ d $ of the previous lemma.
\end{lemma}

\begin{proof}
	The proof is essentially the same as for the previous lemma: the probability that a cluster passes
	between two close points $ z $ and $ z' $ is small (say $ C(z, z') $) for the same reasons. 
	To control the expectation, we can use BK inequality which gives that the probability that
	$ n $ disjoint clusters pass between $ z $ and $ z' $ is smaller that $ C(z, z')^n $.

\begin{figure}[!ht]
\centering
\includegraphics[width=12cm]{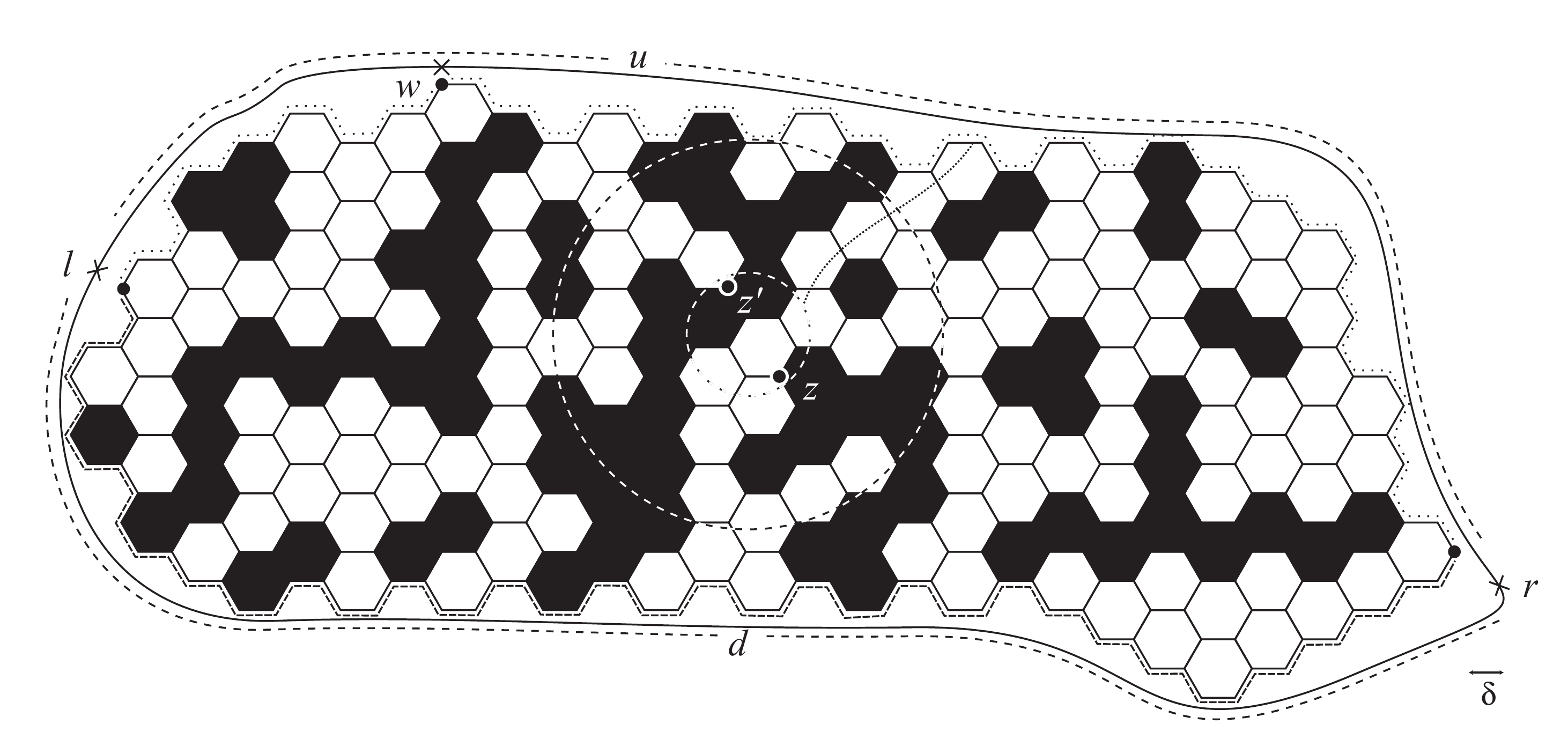}
\caption{The event $ \{ \Nld ( z ) = \Nld ( z' ) + 1 \} $ implies a white connection of a microscopic
circle of diameter $ d ( z, z') $ to a macroscopic circle of diameter $ \beta $ ($ K $ is not fixed precisely 
on this picture).}
\end{figure}

\end{proof}

\begin{proposition}\label{reg}
	The function family $ ( \Hd )_{\delta > 0} $ 
	is precompact with respect to the topology of uniform convergence on every
	compact subset of $ \bar{\Omega} \setminus \{ l, r \}$.
\end{proposition}

\begin{proof}
	We are only interested in letting $ \delta $ tend to $ 0 $ (and otherwise it is anyway trivial).
	So let $ \delta_n $ be a sequence tending to $ 0 $. On each compact subset $ K $ of $ \bar{\Omega} \setminus
	\{ l, r \} $, the functions
	$ \Hld, \Hrd, \Hud, \Hdd $ are bounded and uniformly H\"{o}lder continuous in $ \delta $, 
	so they form equicontinuous families. By Arzel\`{a}-Ascoli's theorem, they form a precompact family. We can
	therefore extract a subsequence $ \delta_k $ of $ \delta_n $ such that $ \Hldk, \Hrdk, \Hudk, \Hddk $ 
	converge uniformly on $ K $. Since $ \bar{\Omega} \setminus \{ l, r \}$ 
	can be written as a countable union of compact subsets, a diagonal extraction gives us the desired result.
\end{proof}

\section{Boundary conditions}

\begin{lemma}\label{hudbound}
	We have the following boundary conditions:
	\begin{eqnarray*}
		\lim_{\delta \to 0} \Hud ( z ) = 0, & & \lim_{\delta \to 0} \Hdd ( z ) = \frac{1}{2},
		\quad \forall z \in \uside \\
		\lim_{\delta \to 0} \Hud ( z ) = \frac{1}{2}, & & \lim_{\delta \to 0} \Hdd ( z ) = 0,
		\quad \forall z \in \dside  \\
		\lim_{\delta \to 0} \Hld (w) & = & \lim_{\delta \to 0} \Hrd (w) = 0
	\end{eqnarray*}
\end{lemma}

\begin{proof}
	By definition and continuity the condition for $ \Hld $ and $ \Hrd $ is obvious.
	
	For the first boundary value, notice that for $ z $ on $ \uside $, the event $ \Qud ( z ) $ implies the connection
	of $ z $ to $ \dside $ (which is at a positive distance from $ z $) by two white paths. By Russo-Seymour-Welsh,
	this probability tends to $ 0 $ as $ \delta \to 0 $, so we are done.
	
	For the second one, first notice that both $ \Qdd (z) $ and its color-negative 
	$ \tilde{\Qdd} (z)$ cannot occur 
	simultaneously for $ z \in u $. 
	By symmetry we obtain that $ \Hdd (z) \leq \frac{1}{2} $. To see that the limit is
	actually $ \frac{1}{2} $, it suffices because of the symmetry to prove that the probability that neither 
	$ \Qdd (z) $ nor $ \tilde{\Qdd} (z) $
	occur tends to $ 0 $ as $ \delta \to 0 $. 
	
	Indeed if $ \Qdd (z) $ does not occur, then either there is no black path separating $ z $ from $ \dside $
	(call this event $ A $) or there is at least one black path separating $ z $ from $ \dside $ but these black
	paths do not touch $ \dside $ (event $ B $).  By self-duality, $ A $ is the event that
	$ z $ is connected to $ \dside $ by a white path. Again by self-duality, the occurence of $ B $ implies that
	$ \tilde{\Qdd} (z) $ occurs: take the lowest black path $ \gamma $ 
	separating $ z $ from $ \dside $ (which does not touch  
	$ \dside $ by definition), so its lower boundary is a white path that touches $ \dside $ (otherwise this white
	path would have a lower boundary which would be a black path and would thus contradict the definition of
	$ \gamma $), which implies that $ \tilde{\Qdd} (z) $ occurs.
	
	So if neither $ \Qdd (z) $ nor $ \tilde{\Qdd} (z) $ happen, $ A $ happens.
	But as seen above, the probability of $ A $ tends to $ 0 $, since the probability of a connection by 
	a white path from $ z $ to $ \dside $ tends to $ 0 $ as $ \delta \to 0 $. 

	The arguments for $ z \in \dside $ are the same as the ones for $ z \in \uside $.
	
	\begin{figure}[!ht]
	\centering
	\includegraphics[width=8cm]{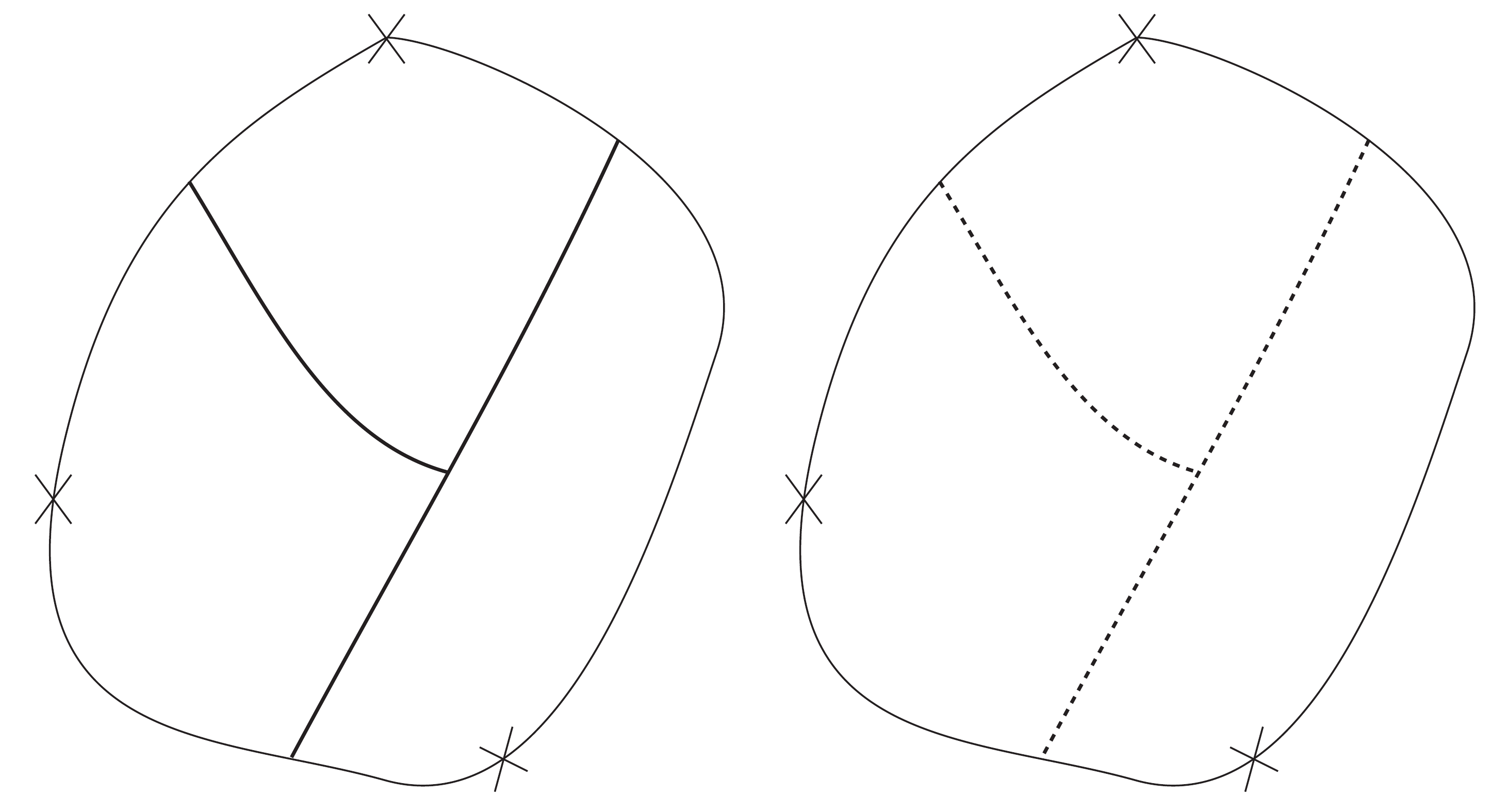}
	\caption{With probability tending to one as $ \delta \to 0 $, exactly one of these two events occurs.}
	\end{figure}
\end{proof}

\section{Analyticity}

We are now interested in showing the analyticity of any subsequential limit of 
$ \Hd = (\Hld + \Hrd) - \hcst (\Hud - \Hdd) $ as
$ \delta \to 0 $ (since by Proposition \ref{reg} the family of functions $ \left( \Hd \right)_{\delta>0} $ is
precompact). The main step consists in proving that for
each $ \delta > 0 $, the function $ \Hd $ is discrete analytic in a sense explained in the next paragraph, which
allows to show that Morera's condition is satisfied.
\subsection{Discrete Cauchy-Riemann equations}

Let us first introduce several notations. 

For an oriented edge $ \onedge = \langle x, y \rangle $ in the interior of $ \Omdel $, let us denote by
\emph{$ \jedge $} and \emph{$ \jsedge $} the edges of $ \Omdel $ obtained by rotating
counterclockwise $ \onedge $  around $ x $ 
by an angle of $ 2 \pi / 3 $ and $ 4 \pi / 3 $ respectively. We will denote by $ \oes $ the \emph{dual edge} of $ \onedge $:
the edge from the center of the hexagon on the right of $ \onedge $ to the center of the hexagon on the left of
$ \onedge $. 

\begin{figure}[!ht]
\centering
\includegraphics[width=6cm]{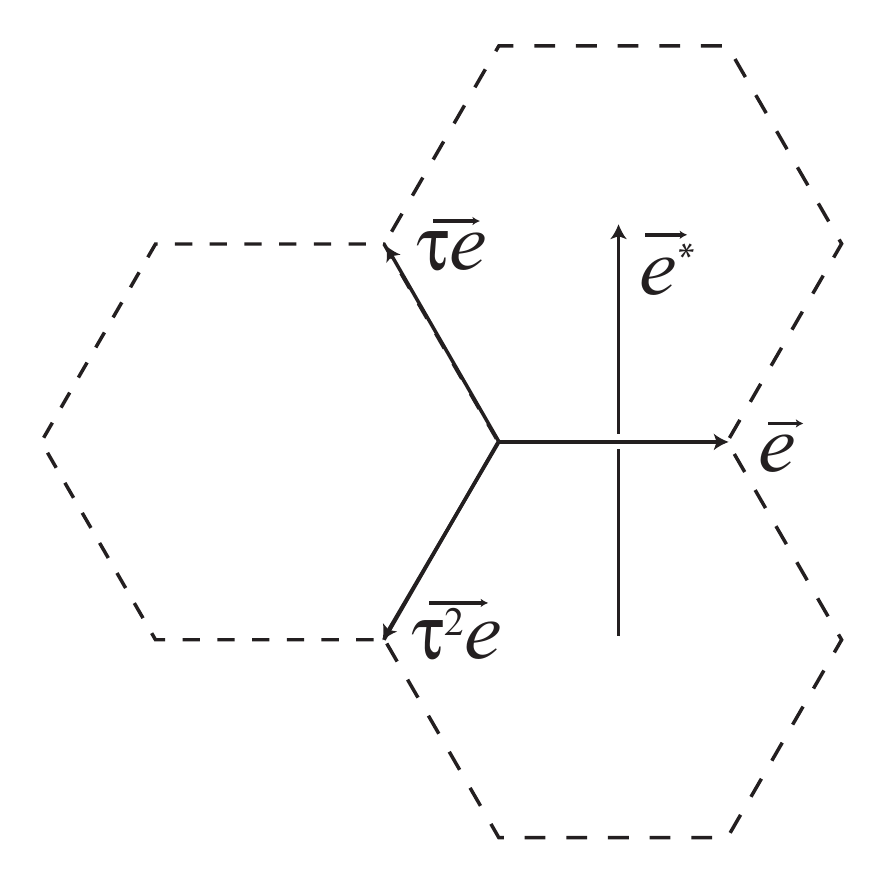}
\caption{Edges notations}
\end{figure}

For a function $ F $ defined on the set of vertices of $ \Omdel $ and an oriented
edge $ \onedge = \langle x, y \rangle $, let us define $ \diffoe F $ as $ F ( y ) - F ( x ) $.

Let $ \diffoepm \Hld $ be $ \Prob [ \Nld ( y ) = \Nld ( x ) \pm 1 ] $ 
By linearity of the expectation it is 
easy to see that $ \diffoe \Hld = \diffoep \Hld - \diffoem \Hld $. 
Let $ \diffoep \Hud $ be $ \Prob [ \Qud ( y ) \setminus \Qud ( x ) ] $ and
$ \diffoem \Hud $ be $ \diffmoep \Hud = \Prob [ \Qud ( x ) \setminus \Qud ( y ) ]$. As before
we have $ \diffoe \Hud = \diffoep \Hud - \diffoem \Hud $. 

For $ \Hrd $ and $ \Hdd $, we define $ \diffoep $ and $ \diffoem $ in the same way as for
$ \Hld $ and $ \Hud $ respectively and also obtain $ \diffoe = \diffoep - \diffoem $. By linearity
it is also defined for $ \Hd $. 

We have the following discrete analyticity result, which already suggests that $ \Hd $ 
is analytic in the limit and is a discrete analogue of the Cauchy-Riemann equations.
\begin{proposition}[Discrete Cauchy-Riemann equations]\label{dcreqs}
	For any $ \delta > 0 $ and any oriented edge $ \onedge $ in the interior of $ \Omdel $, we have the following
	identity:
	\[
		2 \left( \diffoep \Hld - \diffoem \Hrd -  \right)
		= \left( \diffjep - \diffjsep \right) \left( \Hdd - \Hud \right)
	\]
\end{proposition}

\begin{proof}
	Notice that since each configuration (coloring of the hexagons) has equal probability, bijective maps
	are measure-preserving. We will use this fact several times in the proof. Fix $ \delta > 0 $, 
	take as before $ \onedge = \langle x, y \rangle $ and introduce the following notations. 
	
	In what follows, $ \tau . y $ and $ \tau^2 . y $ will be the vertices of $ \Omdel $ such that
	$ \jedge = \langle x, \tau y \rangle $ and $ \jsedge = \langle x, \tau^2 y \rangle $.
	Let $ \LHex $ (respectively $ \RHex $ ; $ \IHex $ ; $ \THex$) 
	be the hexagonal face that is adjacent to $ \onedge $ and 
	$ \jsedge $ (respectively to $ \onedge $ and $ \jedge $ ; to $ \jedge $ and $ \jsedge $ ;
	the hexagon that touches $ y $).
	For a hexagonal face, for instance $ \LHex $, we denote by $ \LHex_w $ the event that this face is connected
	by a white path to $ \dside $, by $ \LHex^b $ the event that it is connected by a black path to $ \uside $,
	by $ \LHex_w^w $ the event that it is connected by (not necessarily disjoint) 
	white paths to both $ \uside $ and $ \dside $, and etc.:
	the connections to $ \uside $ are denoted by superscripts, the connections to $ \dside $ by subscripts.
	Recall that we use the notation $ A \circ B $ for the event that both $ A $ and $ B $ occur on disjoint sites
	(notice that it is well defined for the events we use here).
	
	We now compute the derivative $ \diffoem $ of $ \Hrd $.
	We have that the event $ A := \{ \Nrd (x) = \Nrd (y) + 1 \} $ is the same as $ B:= \Ibb \circ \LHex_w \circ R^w $,
	since it is clear that $ B $ implies $ A $, and by self-duality, if $ B $ does not occur, then $ A $
	does not occur (since otherwise there would be a white path touching the right boundary of the white cluster
	passing between $ y $ and $ x $ and separating $ x $ from $ r $ which would be absurd by definition of 
	the right boundary), so both are equal.

	Notice that on this event, by going from $ y $ to $ x $, either we gain a cluster boundary counting positively
	or we lose a cluster boundary counting negatively.
	
	If $ B $ occurs, then we can define $ \lambda $ as the counterclockwise-most extremal
	white path that joins $ \LHex $ to $ \dside $ (call $ \lambda_d $ its hexagon on $ \dside $) 
	and $ \rho $ as the clockwise-most extremal white path that joins $ \RHex $ to $ \uside $ 
	(call $ \rho_u $ its hexagon on $ \uside $). We can then us a self-duality argument in the interior of rectangle
	$ l, \lambda_d, \THex, \rho_u $ (we consider the topological rectangle delimited by $ \lambda $ (excluded), 
	$ \rho $ (excluded), the arc $ \rho_u l $ (included) and the arc $ l \lambda_d $ (included)):
	$ B $ is the disjoint union of $ C $ and $ D $, where $ C $ is the event that $ B $ happens and that
	there is a white path that joins the arcs $ l \lambda_d $ and $ \THex \rho_u $  
	and $ D $ is the event that
	$ B $ happens and that there is a black path that joins the arcs $ \rho_u l $ and $ \lambda_d \THex $
	(these events occur in the interior of the rectangle).
	So we have $ \Prob [ B ] = \Prob [ C ] + \Prob [ D ] $. But $ C $ is equal to 
	$ \Ibb \circ \LHex_w \circ \Rww $ and we have that $ D $ and $ \Ibb \circ \Lww \circ \RHex^w $ are
	clearly in bijection: it suffices to flip (i.e. invert) the colors inside the rectangle to map one onto the other
	(this is well-defined because the definition of the rectangle does not depend on the colors of the
	hexagons inside), and so the configuration inside is independent of the colors elsewhere.

	But now we have that $ \Ibb \circ \LHex_w \circ \Rww $ and $ \Iww \circ \LHex_b \circ \Rww $ also have the
	same probability. Let $ \iota $ be the clockwise-most extremal black path that joins $ \IHex $ to $ u $, and flip
	the colors in the interior of the part of the graph $ G $ comprised between $ \iota $ and $ \lambda $ that
	contains $ l $ ($ \iota $ and $ \lambda $ excluded). Then flip all the colors of $ \Omdel $. This defines a
	(clearly bijective) map from $ \Ibb \circ \LHex_w \circ \Rww $ to  $ \Iww \circ \LHex_b \circ \Rww $.
	The same color-flipping argument shows that $ \Ibb \circ \Lww \circ \RHex^w $ and $ \Iww \circ \Lww \circ \RHex^b $
	also have the same probability. So we can summarize the discussion above in the following equations, see Figure 10:
	\begin{eqnarray*}
		\diffoem \Hrd & = & \Prob [\Nrd (x) = \Nrd (y) + 1] \\
		& = & \Prob [\Ibb \circ \LHex_w \circ \RHex^w] \\
		& = & \Prob [\Ibb \circ \LHex_w \circ \Rww] + \Prob [\Ibb \circ \Lww \circ \RHex^w] \\
	 	& = & \Prob [ \Iww \circ \LHex_b \circ \Rww ] + \Prob [ \Iww \circ \Lww \circ \RHex^b].
	\end{eqnarray*}

	Using a very similar method (but considering this time a rectangle that contains $ r $ instead of $ l $
	when applying self-duality), one obtains, see Figure 11:
	\begin{eqnarray*}
		\diffoep \Hld & = & \Prob [ \Nld ( y ) = \Nld ( x ) + 1 ] \\
		& = & \Prob [ \Ibb \circ \LHex^w \circ \RHex_w ] \\
		& = & \Prob [ \Ibb \circ \Lww \circ \RHex_w ] + \Prob [ \Ibb \circ \LHex^b \circ \Rww ] \\
		& = & \Prob [ \Iww \circ \Lww \circ \RHex_b ] + \Prob [ \Iww \circ \LHex^b \circ \Rww ]
	\end{eqnarray*}
	
	\begin{figure}[!ht]
	\centering
	\includegraphics[width=12cm]{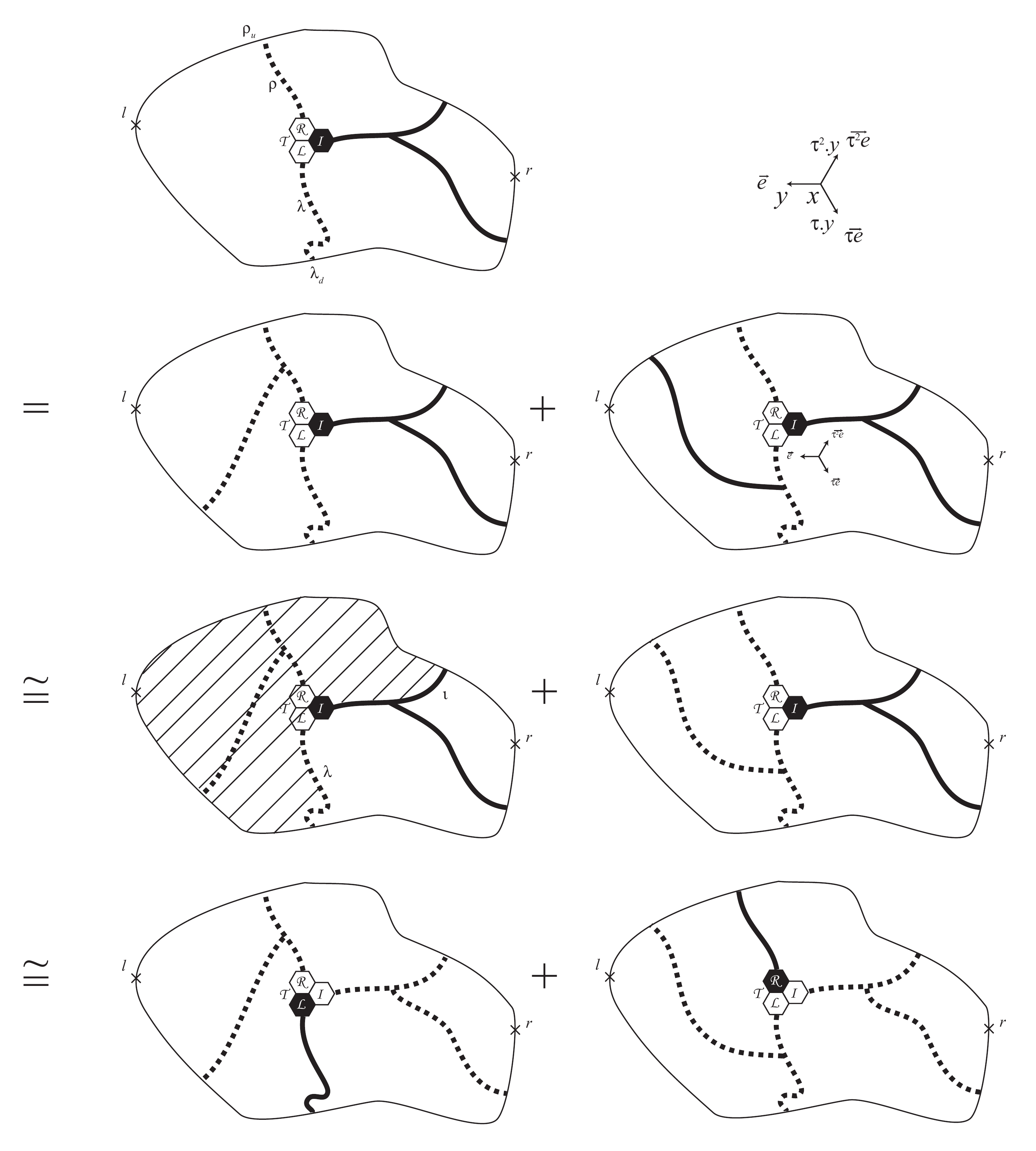}
	\caption{Computation of $ \diffoem \Hrd $. White paths are dashed and black path bold. The stripped
	region is $ G $ (the subgraph where the color are flipped).}
	\end{figure}
	
	\begin{figure}[!ht]
	\centering
	\includegraphics[width=12cm]{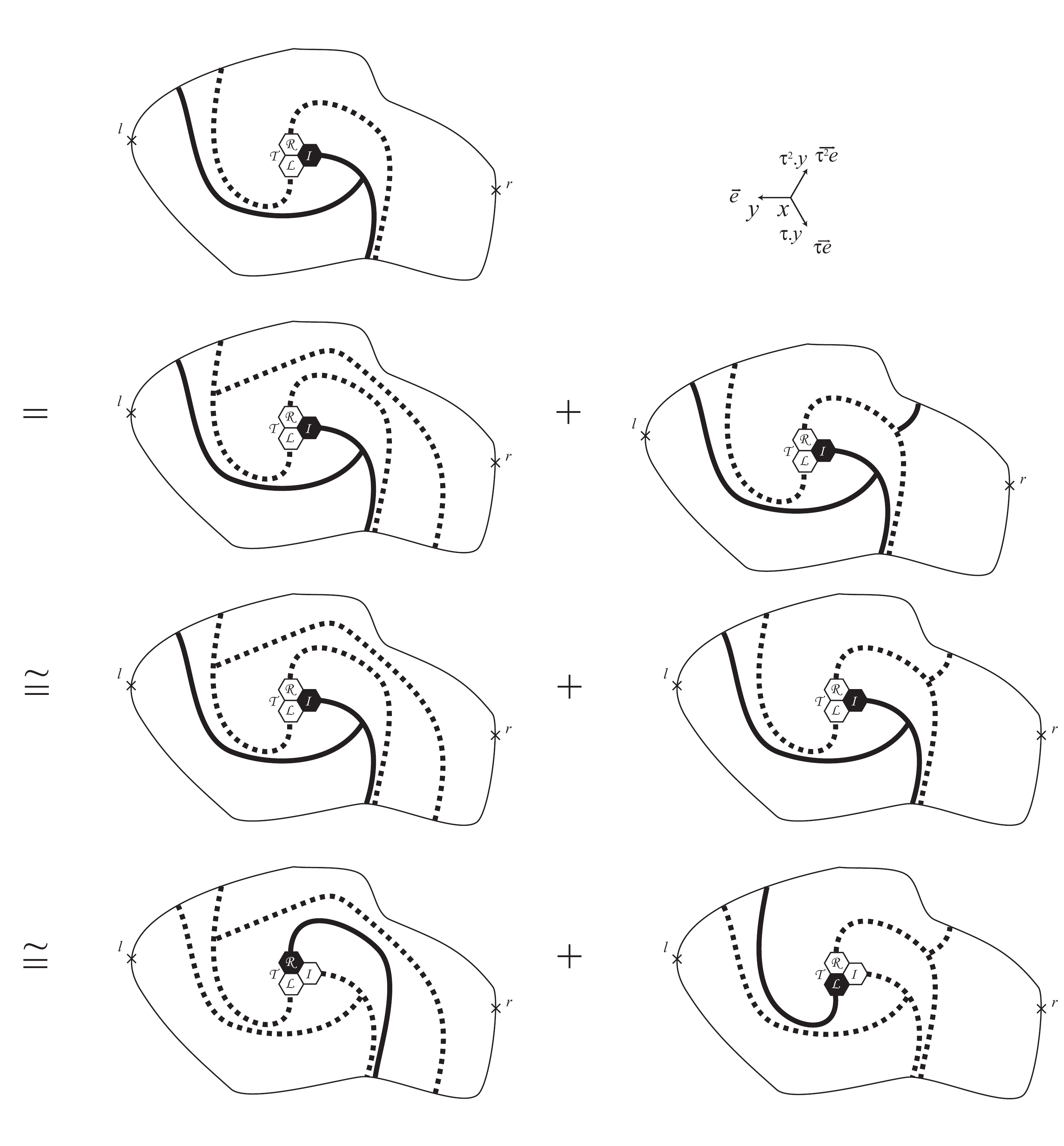}
	\caption{Computation of $ \diffoep \Hld $.}
	\end{figure}

	Let us now compute the derivative $ \diffjep $ of $ \Hud $. By self-duality we have that the event 
	$ X := \Qud ( \tau . y ) \setminus \Qud ( x ) $ is the same as the event that $ \IHex $ and $ \LHex $
	are on a white simple path from $ \dside $ to $ \dside $ which is connected to $ \uside $.
	$ \RHex $ is connected by a black path to $ \uside $ (otherwise there would be a white path separating it
	from $ \uside $ and this path would be connected by a white path to $ \uside $ as well because 
	$ \Qud ( \tau . y ) $ 
	occurs, which would imply that $ \Qud ( x ) $ also occurs). Suppose that $ X $ occurs. Let $ \lambda' $ be 
	the clockwise-most extremal white path that joins $ \LHex $ to $ \dside $ and $ \iota' $ the counterclockwise-most
	extremal white path that joins $ \IHex $ to $ \dside $. Then obviously, exactly one of the
	three following events occurs:
	\begin{enumerate}
		\item $ Y := \Iww \circ \Lww \circ \RHex^b $: 
			there is a white path that joins $ \lambda' $ to $ \uside $ and there is a white path
			that joins $ \iota' $ to $ \uside $.
		\item $ Z := \Iww \circ \LHex_w \circ \Rbb $: there is a white path that joins $ \lambda' $ 
			to $ \uside $ and there is no 
			white path that joins $ \iota' $ to $ \uside $.
		\item $ W := \IHex_w \circ \Lww \circ \Rbb  $: there is no white path that joins 
			$ \lambda' $ to $ \uside $ and there is a
			white path that joins $ \iota' $ to $ \uside $.
	\end{enumerate}
	Using a color-flipping argument we obtain that $ \Prob [Z] = \Prob [\Iww \circ \LHex_b \circ \Rww]$: 
	take the counterclockwise-most black path
	that joins $ \RHex $ to $ \uside $, call it $ \gamma_1 $, the clockwise-most white path that joins $ \LHex $ to 
	$ \dside $, call it $ \gamma_2 $, flip the colors in the interior of the part of $ \Omdel$ delimited by 
	$ \gamma_1 $ and $ \gamma_2 $ that contains $ r $ ($ \gamma_1 $ and $ \gamma_2 $ excluded), then flip the colors
	of the whole graph. This defines a bijection from $ Z $ to $ \Iww \circ \LHex_b \circ \Rww $. Thus we have 
	shown that
	\begin{eqnarray*}
		\diffjep \Hud & = & \Prob[\Iww \circ \Lww \circ \RHex^b] 
		+ \Prob[\Iww \circ \LHex_b \circ \Rww] 
		+ \Prob[\IHex_b \circ \Lww \circ \Rww].
	\end{eqnarray*}

	One obtains similarly, see Figures 13-15 at the end of the section:
	\begin{eqnarray*}
		\diffjep \Hdd & = & \Prob[\Iww \circ \LHex^b \circ \Rww] 
		+ \Prob[\Iww \circ \Lww \circ \RHex_b] 
		+ \Prob[\IHex_b \circ \Lww \circ \Rww] \\
		\diffjsep \Hud & = & \Prob[\Iww \circ \Lww \circ \RHex_b] 
		+ \Prob[\Iww \circ \LHex^b \circ \Rww] 
		+ \Prob[\IHex^b \circ \Lww \circ \Rww] \\
		\diffjsep \Hdd & = & \Prob[\Iww \circ \LHex_b \circ \Rww] 
		+ \Prob[\Iww \circ \Lww \circ \RHex^b] 
		+ \Prob[\IHex^b \circ \Lww \circ \Rww]
	\end{eqnarray*}
	Summing up the identities obtained so far, we obtain the desired result.
\end{proof}

\subsection{Morera's condition}
The last step in order to prove the analyticity of $ h $ is to show that any contour integral of 
the subsequential limit $ h $ vanishes. This
is given by the following proposition (since the convergence is uniform on each compact subset of
$ \Omega $, the integral is equal to the limit of the integrals $ \oint_{\gamma} \Hd (z) \dz $ as
$ \delta \to 0 $).

\begin{proposition}[Morera's condition for $ h $]\label{hmorera}
	Let $ \gamma $ be a simple closed smooth curve in $ \Omega $ oriented counterclockwise. Then we have
	\[
		\oint_{\gamma} \Hd (z) \dz \to 0 \quad \mbox{as $ \delta \to 0 $}
	\]
\end{proposition}
\begin{proof}
	For each sufficiently small $ \delta > 0 $, let $ \gamdel $ be a discretization of $ \gamma $,
	such that $ \gamdel  $ is a simple curve oriented in the same direction consisting in edges
	that follow the orientation of $ \gamdel $ and such that $ \gamdel \to \gamma $ as $ \delta \to 0 $
	in the Hausdorff metric and with a number of edges of order
	$ \delta^{- 1} $. 

	For $ \onedge = \langle x, y \rangle $, let us define $ F ( \onedge ) := \frac{F(x) + F(y)}{2} $ and
	$ \onedge = y - x $ (when appears alone).
	We approximate the integral $ \oint_\gamma \Hd (z) \dz $ by a Riemann sum along 
	$ \gamdel $ defined as $ \sum_{\onedge \in \gamma} \onedge \Hd ( \onedge ) $. 
	
	As $ \delta \to 0 $, one has 
	$ \left| \oint_{\gamma} \Hd (z) \, \mathrm{d}z - \sum_{\onedge \in \gamma} \onedge \Hd ( \onedge ) \right|
	\to 0 $, by precompactness of the family $ \left( \Hd \right)_{\delta>0} $ in the topology of uniform convergence
	on the compact subsets.
	
	We now use the following discrete summation lemma (cf. \cite{Beffara}).
	Define $ \gamint $ as the set of all oriented edges lying in the interior of the part of 
		$ \Omdel $ which is inside $ \gamdel $ and 
	recall that $ \oes $ is the dual edge of $ \onedge $ (seen as a scalar it is equal to 
	$ \sqrt{3} i \onedge $).	
	\begin{lemma}\label{byparts}
		\[
			\sum_{\onedge \in \gamdel} \onedge \Hd ( \onedge ) 
			= \sum_{\onedge  \in \gamint} \oes \diffoep \Hd
			+ o_{\delta \to 0} ( 1 ) 
		\]
	\end{lemma}
	\begin{proof}
		Denote by $ \gamhex $ the set of hexagonal faces of $ \Omdel $ which are inside $ \gamdel $ and
		for such a face $ f $, denote by $ \hexboundary $ the set of its six edges oriented in 
		counterclockwise direction.
		We have that 
		\[
			\sum_{\onedge \in \gamdel} \onedge \Hd ( \onedge )
			= \sum_{f \in \gamhex} \sum_{f \in \hexboundary} \onedge \Hd ( \onedge ),
		\]
		since the terms appearing in edges that are not on $ \gamdel $ appear
		twice (in two faces to which such an edge belongs) 
		with opposite signs and therefore cancel. Denote by 
		$ \langle x_0, x_1 \rangle , \langle x_1, x_2 \rangle, \ldots, \langle x_5, x_0 \rangle $ 
		the six edges of $ \hexboundary $ and take the indices modulo $ 6 $; denote by 
		$ c ( f ) $ the center of a hexagonal face (this term is purely artificial yet). A simple
		calculation shows:
		\[
			\sum_{\onedge \in \hexboundary} \onedge \Hd ( \onedge ) 
			= \sum_{k = 0}^{5} \left( \frac{x_k + x_{k + 1}}{2} - c (f) \right) 
			( H (x_{k + 1}) - H (x_k) ).
		\]
		
		If $ \langle x_k, x_k + 1 \rangle $ does not lie on $ \gamdel $, the term 
		$ \frac{x_k + x_{k + 1}}{2} ( H (x_{k + 1}) - H (x_k) ) $ appears twice with opposite signs
		and cancels, so only the terms with the factor $ c ( f ) $ remain. A term of the form
		$ H (x_{k + 1}) - H (x_k) $ becomes a factor of the difference between two center faces which is
		the edge dual to $ \langle x_k, x_{k + 1} \rangle $. 
		
		On the other hand, we have that the contribution of the boundary terms on $ \gamdel $
		tends to $ 0 $: we have that the number of edges of $ \gamdel $ is of order $ \delta^{- 1} $,
		the term $ \frac{x_k + x_{k + 1}}{2} - c (f) $ is of order $ \delta $ and $ \Hd $ is
		H\"{o}lder on a neighborhood of $ \gamma $.
		
		We obtain that the sum is equal to
		\[
			\sum_{\onedge  \in \mbox{CcwInt} ( \gamdel ) } \oes \diffoe \Hd
			+ o( 1 ), \quad \mbox{as $ \delta \to 0 $}
		\]
		where $ \mbox{CcwInt} $ is the set of the counterclockwise oriented edges of the set of faces
		$ \gamhex $. Taking the sum over the set $ \gamint $ of all oriented edges inside $ \gamdel $,
		using $ \diffoe = \diffoep - \diffmoep $, we obtain
		\[
			\sum_{\onedge  \in \gamint} \oes \diffoep \Hd + o( 1 ), \quad \mbox{as $ \delta \to 0 $}
		\]
		as required.
		
	\end{proof}
	Now it suffices to prove that the sum $ \sum_{\onedge  \in \gamint} \oes \diffoep \Hd $ given by
	the previous lemma is equal to $ 0 $. This is given by the discrete Cauchy-Riemann equations.
	Let us reorder the terms in the sum in the following way:
	\begin{eqnarray*}
		\sum_{\onedge  \in \gamint} \oes \diffoep \Hd
		& = & \sum_{\onedge  \in \gamint} \oes
		\diffoep \left( \Hrd + \Hld - \hcst ( \Hud - \Hdd ) \right) \\
		& = & - \sum_{\onedge  \in \gamint} \oes
		\left( \diffoem \Hrd - \diffoep \Hld + \hcst \diffoep ( \Hud - \Hdd ) \right) \\
		& = & - \sum_{\onedge  \in \gamint} \oes
		\left( \frac{\diffjep \Hud - \diffjep \Hdd - \diffjsep \Hud + \diffjsem \Hdd}{2} \right. \\
		& & \left. + \hcst \diffoep ( \Hud - \Hdd ) \right),
	\end{eqnarray*}
	where last equality is obtained using the discrete Cauchy-Riemann equations. Reordering one last time the
	sum (using the changes of variables $ \jes \to \oes $ and $ \jses \to \oes $ in the first and
	the second parts of the sum respectively), we obtain
	\[
		-\frac{1}{2} \sum_{\onedge \in \gamint} \left(\sqrt{3} i \oes + (\jses - \jes)  \right) \diffoep ( \Hud - \Hdd ),
	\]
	which is equal to $ 0 $, since $ \left(\sqrt{3} i \oes + (\jses - \jes)  \right) = 0 $ by the geometry
	of the lattice (and this is in fact the only 
	step in our proof where the actual \emph{embedding} of the lattice is crucial).
	
\end{proof}

\begin{figure}[!htp]
\centering
\includegraphics[width=12cm]{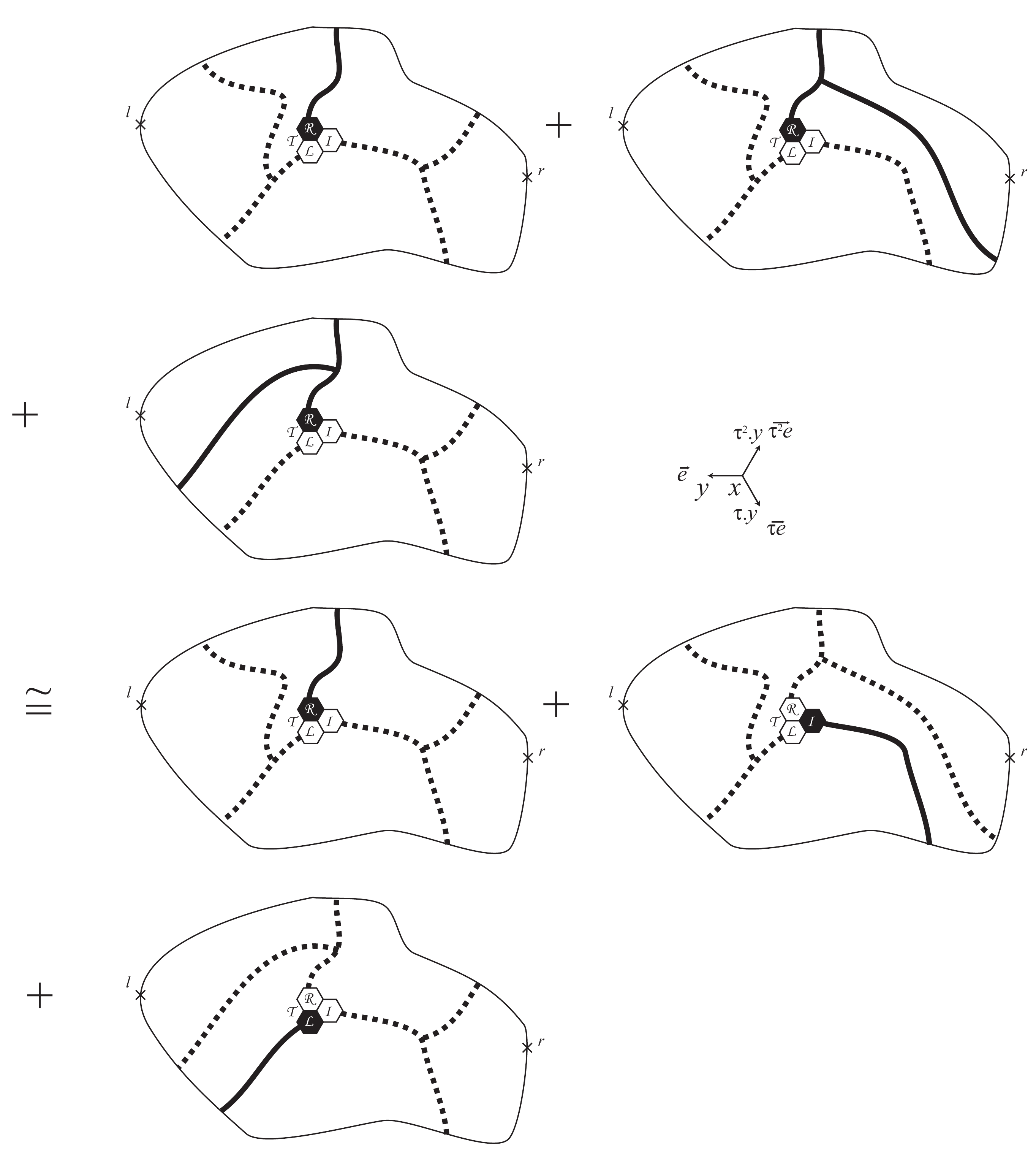}
\caption{Computation of $ \diffjep \Hud $.}
\end{figure}

\begin{figure}[!htp]
\centering
\includegraphics[width=12cm]{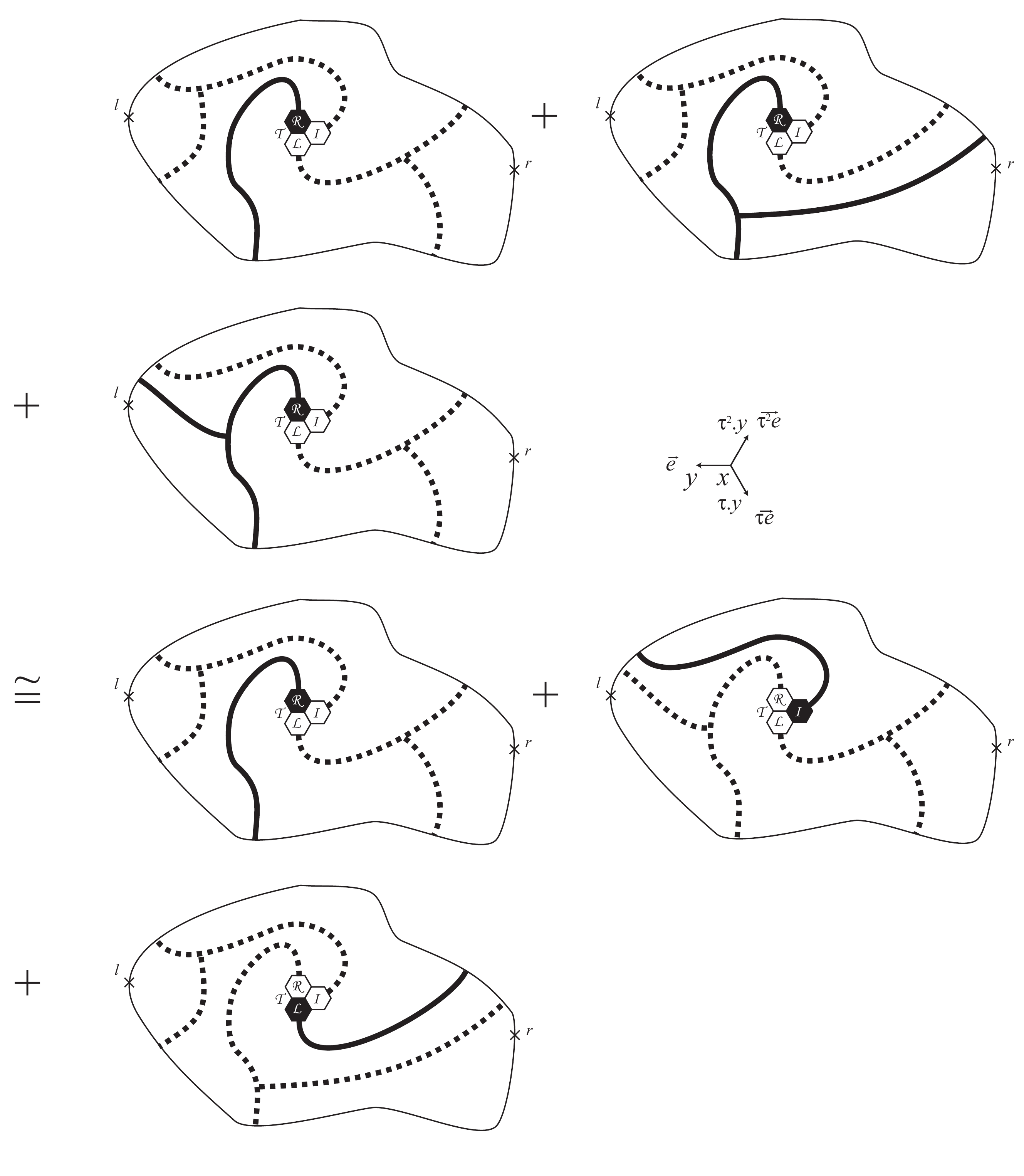}
\caption{Computation of $ \diffjep \Hdd $.}
\end{figure}

\begin{figure}[!htp]
\centering
\includegraphics[width=12cm]{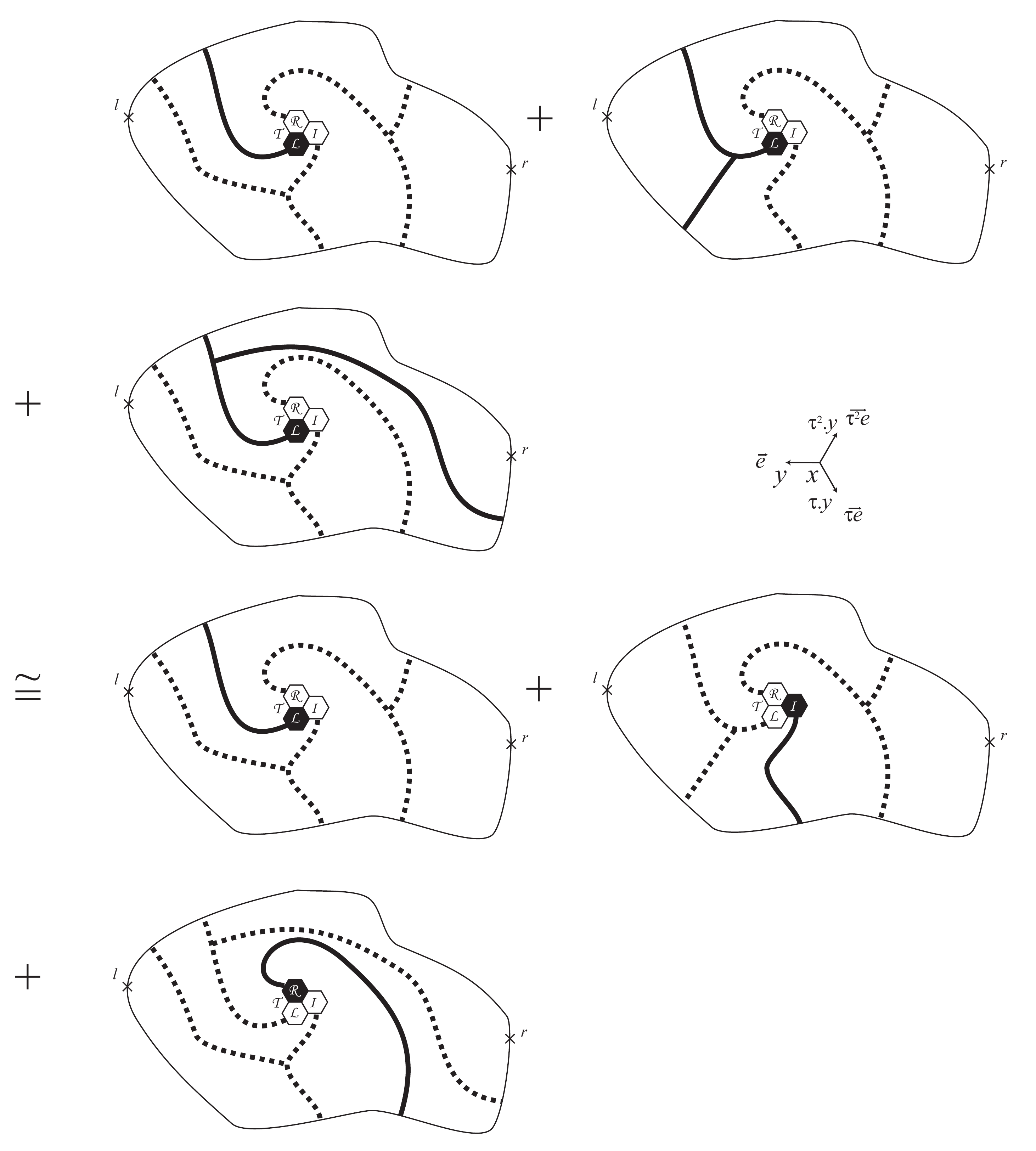}
\caption{Computation of $ \diffjsep \Hud $.}
\end{figure}

\begin{figure}[!ht]
\centering
\includegraphics[width=12cm]{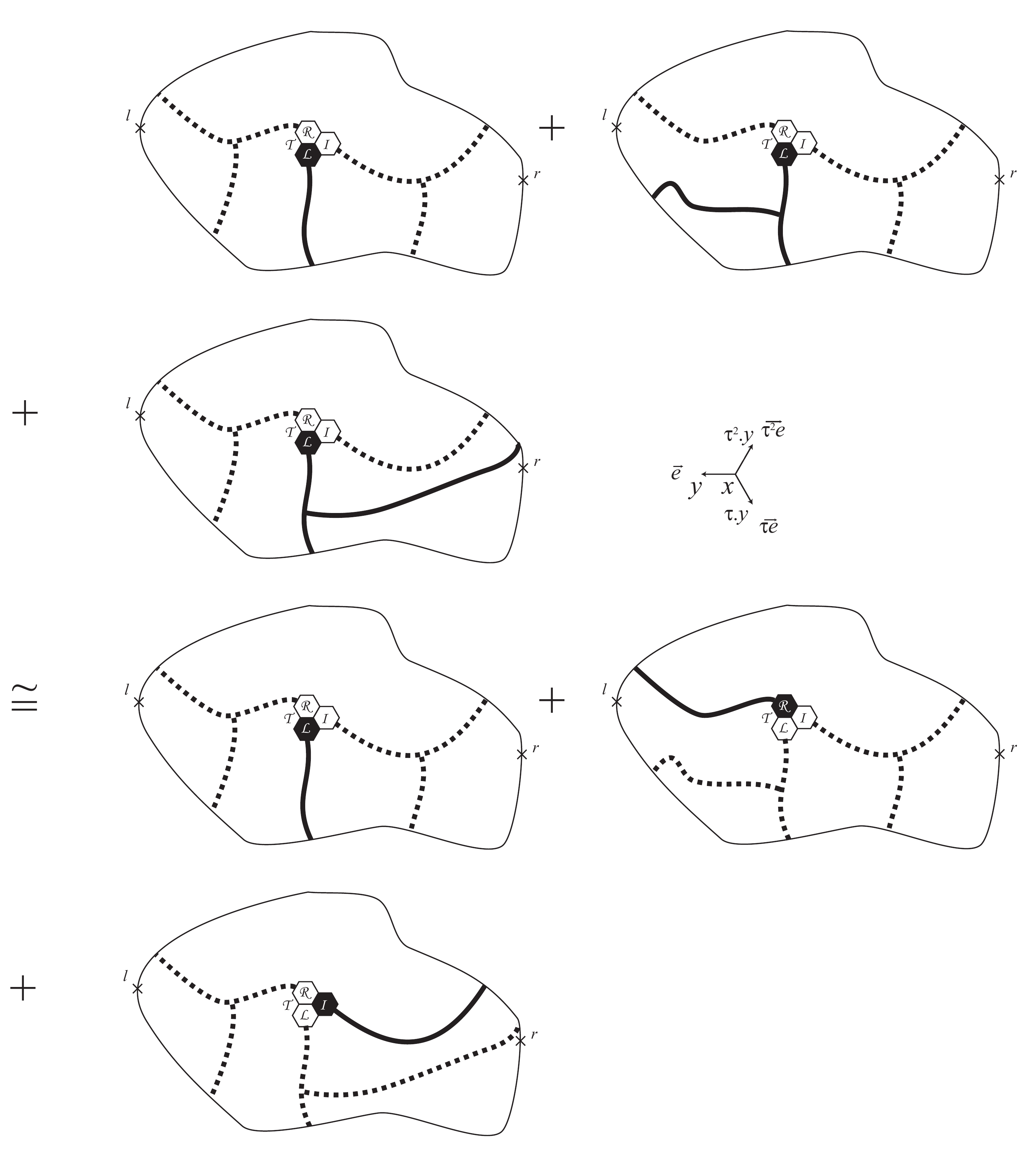}
\caption{Computation of $ \diffjsep \Hdd $.}
\end{figure}


\newpage

\begin{thebibliography}{9}

\bibitem{Beffara}
	V. Beffara, \emph{Cardy's formula on the triangular lattice, the easy way},
in \emph{Universality and renormalization}, 39--45, Fields Inst. Commun., 50, Amer. Math. Soc., Providence, RI, 2007.
\bibitem{Bollobas}
	B. Bollobas, O. Riordan, \emph{Percolation}, Cambridge University Press, 2006.
\bibitem{Cardy}
	J. Cardy, \emph{Critical percolation in finite geometries}, J. Phys. A, 25:L201--206, 1992.
\bibitem{Cardy-expected}
	J. Cardy, \emph{Conformal Invariance and Percolation}, arXiv:math-ph/0103018.
\bibitem{Camia}
	F. Camia, C. Newman, \emph{Critical percolation exploration path and SLE 6: a proof of convergence},
	Probability Theory and Related Fields, 139, (3), 473--519, 2007.
\bibitem{Dubedat}
	J. Dub\'{e}dat, \emph{Excursion decompositions for SLE and Watts' crossing formula}, 
	Probab. Theory Related Fields, 2006, no. 3, 453--488.
\bibitem{Grimmett}
	G. Grimmett, \emph{Percolation}, Springer-Verlag, 1999.
\bibitem{Lawler}
	G. F. Lawler, \emph{Conformally Invariant Processes in the Plane}, volume 114. Mathematical
	Surveys and Monographs, 2005.
\bibitem{Simmons}
	J. H. Simmons, P. Kleban, R. M. Ziff, \emph{Percolation crossing formulae and conformal field theory},
J. Phys. A 40 (2007), no. 31, F771--F784. 
\bibitem{Smirnov}
	S. Smirnov, \emph{Critical percolation in the plane: Conformal invariance, Cardy's formula, scaling limits},
	C. R. Acad. Sci. Paris Sr. I Math. 333 (2001), 239--244.
\bibitem{Smirnov-preprint}
	S. Smirnov, \emph{Critical percolation in the plane},
preprint, 2001.
\bibitem{Smirnov-icmp}
	S. Smirnov, \emph{Critical percolation and conformal invariance},
 XIVth International Congress on Mathematical Physics (Lissbon, 2003), 
 99--112, World Sci. Publ., Hackensack, NJ. 
\bibitem{Smirnov2}
	S. Smirnov, \emph{Towards conformal invariance of 2D lattice models}, Sanz-Sole, Marta (ed.) et al., Proceedings
of the international congress of mathematicians (ICM), Madrid, Spain, August 22-30, 2006. Volume II: Invited
lectures, 1421--1451. Zurich: European Mathematical Society (EMS), 2006.
\bibitem{Smirnov-Werner}
	S. Smirnov, W. Werner, 
\emph{Critical exponents for two-dimensional percolation},
Math. Research Letters 8 (2001), no. 5--6, 729--744, 2001. 
\bibitem{Watts}
	G.M.T. Watts. \emph{A crossing probability for critical percolation in two dimensions}, J. Phys.
	A29:L363, 1996.

\end{thebibliography}
\end{document}